\documentclass[a4paper,10pt,leqno,twoside]{tobart}

\usepackage[english]{babel} \usepackage{inputenc, amsmath, amssymb, latexsym,
  epsfig, graphicx, rotating, fancyhdr, amsthm, pifont, empheq}

\usepackage{color}

\newcommand{\ld}{\ensuremath{,\ldots,}}
\newcommand{\ssq}{\ensuremath{\subseteq}}
\newcommand{\smin}{\ensuremath{\setminus}}
\newcommand{\eps}{\ensuremath{\varepsilon}}
\newcommand{\wh}{\ensuremath{\widehat}}

\newcommand{\inte}{\ensuremath{\mathrm{int}}}

\newcommand{\diam}{\ensuremath{\mathrm{diam}}}

\newcommand{\kreis}{\ensuremath{\mathbb{T}^{1}}}
\newcommand{\ntorus}[1][2]{\ensuremath{\mathbb{T}^{#1}}}

\newcommand{\sltr}{\ensuremath{\textrm{SL}(2,\mathbb{R})}}
\newcommand{\torus}{\ensuremath{\mathbb{T}^2}}

\newcommand{\alphlist}{\begin{list}{(\alph{enumi})}{\usecounter{enumi}\setlength{\parsep}{2pt}
      \setlength{\itemsep}{1pt} \setlength{\topsep}{5pt}
      \setlength{\partopsep}{3pt}}}
\newcommand{\arablist}{\begin{list}{(\arabic{enumi})}{\usecounter{enumi}\setlength{\parsep}{2pt}
          \setlength{\itemsep}{1pt} \setlength{\topsep}{5pt}
          \setlength{\partopsep}{3pt}}}
\newcommand{\romanlist}{\begin{list}{(\roman{enumi})}{\usecounter{enumi}\setlength{\parsep}{2pt}
              \setlength{\itemsep}{1pt} \setlength{\topsep}{5pt}
              \setlength{\partopsep}{3pt}}}
\newcommand{\Romanlist}{\begin{list}{(\Roman{enumi})}{\usecounter{enumi}\setlength{\parsep}{2pt}
              \setlength{\itemsep}{1pt} \setlength{\topsep}{5pt}
              \setlength{\partopsep}{3pt}}}
\newcommand{\bulletlist}{\begin{list}{$\bullet$}{\setlength{\parsep}{2pt}
                \setlength{\itemsep}{1pt} \setlength{\topsep}{5pt}
                \setlength{\partopsep}{3pt}\setlength{\leftmargin}{15pt}}} 
\newcommand{\Alphlist}{\begin{list}{(\Alph{enumi})}{\usecounter{enumi}\setlength{\parsep}{2pt}
      \setlength{\itemsep}{1pt} \setlength{\topsep}{5pt}
      \setlength{\partopsep}{3pt}}}
 \newcommand{\listend}{\end{list}}

\newcommand{\T}{\ensuremath{\mathbb{T}}}

\newcommand{\N}{\ensuremath{\mathbb{N}}} 
\newcommand{\R}{\ensuremath{\mathbb{R}}}
\newcommand{\Z}{\ensuremath{\mathbb{Z}}}
\newcommand{\Q}{\ensuremath{\mathbb{Q}}}

\newcommand{\cA}{\mathcal{A}}
\newcommand{\cB}{\mathcal{B}}
\newcommand{\cC}{\mathcal{C}}

\newcommand{\cF}{\mathcal{F}}

\newcommand{\cP}{\mathcal{P}}

\newcommand{\cS}{\mathcal{S}}
\newcommand{\cT}{\mathcal{T}}

\newcommand{\nLim}{\ensuremath{\lim_{n\rightarrow\infty}}}

\newcommand{\viertel}{\ensuremath{\frac{1}{4}}}

\newtheoremstyle{tobthm}{3pt}{3pt}{\itshape}{0pt}{\bfseries}{.}{0.5eM}{}
\theoremstyle{tobthm}

\newtheorem{definition}{Definition}[section]
\newtheorem{thm}[definition]{Theorem}

\newtheorem{lem}[definition]{Lemma}
\newtheorem{Lemma}[definition]{Lemma}
\newtheorem{cor}[definition]{Corollary}
\newtheorem{prop}[definition]{Proposition}

\newtheorem{claim}[definition]{Claim}

\newtheoremstyle{tobrem}{3pt}{3pt}{\normalfont}{0pt}{\bfseries}{.}{0.5em}{}
\theoremstyle{tobrem}

\newtheorem{rem}[definition]{Remark}

\newcommand{\vertset}{\ensuremath{\operatorname{Vert}}}
\newcommand{\VP}{\ensuremath{\operatorname{VP}}}
\newcommand{\CV}{\ensuremath{\operatorname{CV}}}
\newcommand{\RCV}{\ensuremath{\operatorname{RCV}}}
\newcommand{\LS}{\ensuremath{\operatorname{LS}}}

\numberwithin{equation}{section} \numberwithin{figure}{section}

\title{\large \bf Genericity of mode-locking for quasiperiodically forced circle
  maps} \author{J.~Wang\thanks{Department of Applied Mathematics, Nanjing
    University of Science and Technology, 210094 Nanjing, China. Email: {\tt
      jingwang018@gmail.com}}, Q.~Zhou\thanks{Department of Mathematics, Nanjing
    University, 210093 Nanjing, China. Email: {\tt qizhou@nju.edu.cn}}, and
  T.~J\"ager\thanks{Friedrich-Schiller-University Jena, Institute of
    Mathematics, 07743 Jena, Germany. Email: {\tt tobias.jaeger@uni-jena.de}}}

\begin{document}

\setlength{\abovedisplayskip}{1.0ex}
\setlength{\abovedisplayshortskip}{0.8ex}

\setlength{\belowdisplayskip}{1.0ex}
\setlength{\belowdisplayshortskip}{0.8ex}

\maketitle

\abstract{We show that a generic quasiperiodically forced circle homeomorphism
  is mode-locked: the rotation number in the fibres is rationally related to the
  rotation number in the base and it is stable under small perturbations of the
  system. As a consequence, this implies that for a generic parameter family of
  quasiperiodically forced circle homeomorphisms satisfying a twist condition,
  the graph of the rotation number as a function of the parameter is a devil's
  staircase. }

\section{Introduction} \label{Introduction}

One of the cornerstones of low-dimensional dynamics is Poincar\'e's rotation
theory on the circle. Given an orientation-preserving circle homeomorphism $f :
\kreis \to \kreis$, where $\kreis = \R/\Z$, and a lift $F:\R\to\R$, the rotation
number of $F$ is defined as $\rho(F)=\nLim (F^n(x)-x)/n$. This limit always
exists and is independent of $x$
(e.g.~\cite{demelo/vanstrien:1993,katok/hasselblatt:1997}). In this context,
{\em mode-locking} means that the rotation number of $f$ is constant under small
$\mathcal{C}^0$-perturbations. This can only happen if $\rho(f)$ is rational and
occurs, for example, if $f$ has a stable fixed or periodic point. More
generally, $f$ is mode-locked if and only if there exists a non-empty interval
$I\subsetneq\kreis$ and some $p\in\N$ such that $f^p$ maps $I$ strictly inside
itself (that is $f^p\left(\overline{I}\right)\ssq \inte(I)$).  This is true for an open and
dense set of orientation-preserving circle homeomorphisms, so that in this
setting mode-locking is generic in a topological sense. Further, one may
consider parameter families of circle homeomorphisms $(f_\tau)_{\tau\in\kreis}$
which satisfy a twist condition, that is, for all $x\in\kreis$ the mapping
$\tau\mapsto f_{\tau}(x)$ is strictly order-preserving (that is, an
orientation-preserving circle homeomorphism as well). Then as a consequence of
the above, one obtains that in the topologically generic case the map
$\tau\mapsto\rho(f_\tau)$ is locally constant on an open and dense set of
parameters, so that its graph is a {\em devil's staircase}.

Our aim is to study the analogous phenomenon in the context of {\em
  quasiperiodically forced (qpf) circle homeomorphisms}. These are
homeomorphisms of the two-dimensional torus, homotopic to the identity, with skew product
structure
\begin{equation} \label{e.qpf-system} f:\T^2\to\T^2 \quad , \quad
  f(\theta,x)=(\theta+\omega,f_\theta(x)) \ .
\end{equation}
Usually $\omega\in\kreis$ is supposed to be irrational -- hence quasiperiodic
forcing. However, for our purposes it will be convenient to include the rational
case and to speak of periodically forced circle homeomorphisms in this case. The
circle homeomorphisms $f_\theta:\kreis\to\kreis$ are called {\em fibre
  maps}.

There is one particular case that has to be mentioned in this context, which are
quasiperiodic $\sltr$-cocycles. Via their projective action, these induce maps
of the form (\ref{e.qpf-system}) (e.g.\ \cite{haro/puig:2006}). For the special
case of Schr\"odinger cocycles with potential
$V(\theta)=\lambda\cos(2\pi\theta)$, the existence of a devil's staircase for
the respective parameter families is equivalent to Cantor spectrum of the so
called almost-Mathieu operator, which has been established in
\cite{bellissard/simon:1982,CEY, puig:2004,avila/jitomirskaya:2005}.  For more
general potentials $V$, Simon conjectured that the corresponding one-dimensional
quasiperiodic Schr\"odinger operator generically has Cantor spectrum
\cite{simon}. We recall that one dimensional quasiperiodic Schr\"odinger
operator is
\begin{equation}\label{equ-schrodinger}
(H u)_n= u_{n+1}+u_{n-1} + V(\theta+n\omega ) u_n
\end{equation}
with frequency $\omega\in\T^d (d\in\N),$ and potential function $V\in
C^k(\T^d,\R)$, $(k\in \N, \infty, \omega)$. A number of results have been
obtained on this problem, with the most studied case being the analytic one.
Eliasson proved that for any fixed Diophantine frequency $ \omega \in \T^d$, the
quasiperiodic Schr\"{o}dinger operator with generic small analytic potential has
Cantor spectrum \cite{eliasson:1992} . Then, still in the analytic topology with
$d=1$, Goldstein and Schlag proved that for generic $\omega \in \R\backslash
\Q$, in the region of positive Lyapunov exponents (that is, for $V$ sufficiently
large), the spectrum is again a Cantor set \cite{GS11}.  In the
$\cC^0$-topology, Avila-Bochi-Damanik proved Cantor spectrum for any fixed
totally irrational vector $ \omega \in\T^d$ and generic $V \in C^0(\T^d,\R)$ \cite{ABD,aviladense}.
In the case of $C^k$-topology ($1\leq k\leq \infty$ or even in analytic
category), Simon's conjecture is true for $d=1$: one can prove that for generic $
\omega \in \R\backslash \Q$ and generic $v \in C^k(\T,\R)$, the spectrum of
(\ref{equ-schrodinger}) is a Cantor set.\footnote{ To see this, one can first
  perturb the potential such that its related cocycle has positive Lyapunov
  exponent \cite{aviladense}, and then perturb the frequency to obtain Cantor
  spectrum \cite{GS11}.  This proof has been pointed out to the second author by
  Avila.} Results on the continuous-time case can be found in
\cite{eliasson:1992,FJP}.

Yet, for all these results the respective proofs depend strongly on the
particular linear structure of \sltr-cocycles and can usually not be carried
over directly to a more general nonlinear setting.  A genuinely non-linear
approach that is based on multiscale analysis and parameter exclusion in the
spirit of Benedicks and Carleson \cite{benedicks/carleson:1991} has been
developed in \cite{bjerkloev:2005a,jaeger:2009a,Jaeger2012SNADiophantine} and
applied to the problem of mode-locking both in the linear and the non-linear
case
\cite{WangZhang2013Schroedinger,WangZhang2014CantorSpectrum,JaegerWang2015ModeLocking}.
This allows to provide effective criteria for the existence of a devil's
staircase in explicit parameter families, although up to date only partial
results are available for non-linear systems
\cite{JaegerWang2015ModeLocking}. In particular, there exist no explicit
examples of non-linear systems for which a complete devil's staircase has been
established.\medskip

Here, we take a different route and study the problem from the viewpoint of
topological genericity, analogous to the results on generic Schr\"odinger
operators mentioned above. We denote the class of maps of the form
(\ref{e.qpf-system}) by $\mathcal{F}$ and equip them with the metric
$d(f,g)=\max\left\{d_{\mathcal{C}^0}(f,g),d_{\mathcal{C}^0}(f^{-1},g^{-1})\right\}$,
where $d_{\cC^0}$ denotes the canonical $\cC^0$-distance. Note that thus
$\mathcal{F}$ becomes a complete metric space, and hence a Baire space.  The
{\em fibred rotation number} is given by
\begin{equation}
  \label{e.fibred_rotnum}
  \rho(F,\theta) \ = \ \nLim \left(F^n_\theta(x)-x\right)/n \ ,
\end{equation}
where $F:\kreis\times \R\to\kreis\times\R$ is a lift of $f$ and we write
$F^n_\theta=F_{\theta+(n-1)\omega}\circ \ldots \circ F_\theta$ for the fibre maps of
the iterates. Again, the limit in (\ref{e.fibred_rotnum}) always exists and is
independent of $x\in\R$. Moreover, if $\omega$ is irrational, then it is also
independent of $\theta$ (and will simply be denoted by $\rho(F)$) and the
convergence is uniform on $\T^2$ \cite{johnson/moser:1982,herman:1983}.  The
mechanism for mode-locking in qpf circle homeomorphisms in this situation is
essentially the same as in the unforced case.
\begin{thm}[\cite{bjerkloev/jaeger:2009}] A qpf circle homeomorphism $f$ is
  mode-locked if and only if there exists a topological annulus $\cA\ssq \T^2$
  and some $p\in\N$ such that $f^p$ maps $\cA$ strictly inside itself.
\end{thm}
\begin{rem}
  It can be shown that the annulus $\cA$ is necessarily non-trivial in homotopy,
  so it `wraps aroud the torus' in some direction. Moreover, in the terminology
  of \cite{bjerkloev/jaeger:2009}, the intersection $\bigcap_{n\in\N}f^{np}(\cA)$
  is a $p$-periodic strip, and its existence forces the rotation numbers
  $\omega$ and $\rho(f)$ to be rationally related (that is, the equation
  $r\omega+s\rho(f)+t=0$ has a non-zero integer solution $(\hat r,\hat s,\hat
  t)$). We refer to \cite{bjerkloev/jaeger:2009} for details.
\end{rem}

One advantage of this equivalent characterisation of mode-locking is that it can
directly be extended to periodically forced circle homeomorphisms -- where
mode-locking is more difficult to define otherwise due to the independence of
the fibres over disjoint orbits -- and that it is obviously an open
condition. It is thus perfectly adapted to address questions of genericity. We
therefore make the following definition.
\begin{definition} \label{d.mode-locking}
  We say $f\in\cF$ is {\em mode-locked}, if there exist a topological annulus
  $\cA\ssq\torus$ and an integer $p$ such that $f^p$ maps $\cA$ strictly inside
  itself.
\end{definition}

 Our main result is the following.
\begin{thm} \label{t.main-homeomorphisms} There exists an open and dense set
  $\cF^\textrm{ML}\ssq \cF$ such that all $f\in\cF^\textrm{ML}$ are mode-locked.
\end{thm}

 To the best knowledge of the authors, this result establishes the first
 nonlinear analogue to the Cantor spectrum results discussed above. It is
 important to note here that we work with the whole space $\cF$, which allows to
 change the rotation number and use rational approximations for the
 proof. However, this also means that we obtain a statement mainly about
 Liouvillean rotation numbers in the base (which are topologically generic), and
 nothing can be said in this way about the Diophantine case (which is
 topologically meager). This becomes evident through the following statement.
 Given $\omega\in\kreis$, we let $\cF_\omega=\{f\in\cF\mid f\textrm{ has
   rotation number } \omega \textrm{ in the base } \}$ and equip this space with
 the induced $\cC^0$-topology. Recall that a subset of a Baire space is called
 {\em residual} if it is a countable intersection of open and dense sets.
\begin{cor}\label{c.main-homeomorphisms}
  There exists a residual set $\Omega\ssq\kreis$ such that for all
  $\omega\in\Omega$ there is an open and dense set $\cF_\omega^\textrm{ML}$ such
  that all $f\in\cF^\textrm{ML}_\omega$ are mode-locked.
\end{cor}
As a residual set, $\Omega$ will contain a residual subset of Liouvillean
rotation numbers (since intersections of residual sets are residual), but does
not necessarily contain any Diophantine numbers.

As in the one-dimensional case, we then turn to consider parameter families
$(f_\tau)_{\tau\in\kreis}$ with a twist condition, by which we mean that for all
$(\theta,x)\in\T^2$ the mapping $\tau\mapsto f_{\tau,\theta}(x)$ is strictly
orientation-preserving. We denote the space of such families with fixed rotation
number $\omega$ in the base by $\cP_\omega$, let $\cP=\bigcup_{\omega\in\kreis}
\cP_\omega$ and equip both sets with the respective $\cC^0$-topology.
\begin{thm}\label{t.main-families}
  \begin{itemize}
  \item[(a)] There exists a residual set $\cP^\textrm{DS}\ssq\cP$ such that for
    all $(f_\tau)_{\tau\in\kreis} \in\cP^\textrm{DS}$ the graph of the function
    $\tau\mapsto\rho(f_\tau)$ is a devil's staircase.
  \item[(b)]There exists a residual set $\Omega\ssq\kreis$ with the property
    that for all $\omega\in\Omega$ there is a residual set
    $\cP_\omega^\textrm{DS}$ such that for all
    $(f_\tau)_{\tau\in\kreis}\in\cP^\textrm{DS}_\omega$ the graph of the
    function $\tau\mapsto\rho(f_\tau)$ is a devil's staircase.
  \end{itemize}

\end{thm}

\noindent {\bf Acknowledgments.} We would like to thank the Alexander von Humboldt-Foundation,
who has made this collaboration possible by supporting JW with a postdoc
fellowship during her stay at the TU Dresden and the FSU Jena. TJ acknowledges
support by the German Research Council (Emmy Noether grant Ja 1721/2-1 and
Heisenberg grant Oe 538/6-1) and the EU-funded CRITICS network.

\section{Notation and preliminaries} \label{Preliminaries}

Given $x,y\in\kreis$, we denote by $[x,y]$ the positively oriented closed
interval from $x$ to $y$, and use similar notation for open and half-open
intervals. Maps defined on tori will usually be denoted by small letters (greek
or latin) and their lifts by the respective captital letters. If
$f:\T^2\to\torus$ is a qpf circle homeomorphism, then we usually consider lifts
$F:\T^1\times\R\to\kreis\times \R$ to the open annulus (and only occasionally
lifts to the plane).

It will be convenient to introduce the following alternative notation for qpf
circle homeomorphisms. Suppose $\varphi:\torus\to\kreis$ is a continuous
function such that for all $\theta\in\kreis$ the mapping $x\mapsto
x+\varphi(\theta,x)$ defines an orientation-preserving circle
homeomorphism. When $\varphi$ is differentiable, the latter means that
$\partial_x\varphi(\theta,x)\geq-1$ for all $x$ and the set
$\{x\in\kreis\mid \partial_x\varphi(\theta,x)=-1\}$ is nowhere dense. Here, we
use $\partial_\vartheta$ to denote the derivative with respect to the variable
$\vartheta$.

Given $\omega\in\kreis$, we can then define a map $(\omega,\varphi)$ of the form
(\ref{e.qpf-system}) by
\begin{equation}
  (\omega,\varphi) \ : \ (\theta,x)\mapsto (\theta+\omega,x+\varphi(\theta,x)) \ .
\end{equation}
We denote the space of functions $\varphi:\torus\to\kreis$ with the above properties
by $\Xi$. Note that $\Xi$ is `convex' in the following sense: If
$\varphi,\hat\varphi\in\Xi$ and $\Phi,\hat\Phi:\torus\to\R$ are lifts, then any convex
combination of $\Phi$ and $\hat \Phi$ is again a lift of a function in
$\Xi$. Thus, it is always possible to interpolate between two elements of $\Xi$
in a continuous way. In general, the path obtained by this interpolation depends
on the choice of the lifts $\Phi,\hat\Phi$, which are only well-defined modulo
additive integer constants. However, we will often have that $\hat\varphi$ is a
small $\cC^0$-perturbation of $\varphi$, that is, $\eps=d_{\cC^0}(\varphi,\hat\varphi)\ll
1$. In this case, there is a natural way to interpolate between $\varphi$ to
$\hat\varphi$ within the $\eps$-neighbourhood of the two functions, which is simply
to chose lifts $\Phi,\hat\Phi$ with
$d_{\cC^0}(\Phi,\hat\Phi)=d_{\cC^0}(\varphi,\hat\varphi)=\eps$ for the convex combinations. In
this situation, we will slightly abuse notation and write $t\varphi+(1-t)\hat\varphi$
without further explanation.

We will further need the following notation and conventions. Given a qpf circle
homeomorphism $f$ with rotation number $\omega$ on the base, we let
$\varphi_f(\theta,x)=f_\theta(x)-x$, so that $(\omega,\varphi_f)=f$. Further,
given $\varphi\in\Xi$, we let $\varphi^{-1}=\varphi_{(\omega,\varphi)^{-1}}$ denote
the respective translation function of the inverse of $(\omega,\varphi)$, so
that $(-\omega,\varphi^{-1})=(\omega,\varphi)^{-1}$. More generally, we denote
by $\varphi^q$ the translation function associated to $f^q$, that is, $f^q=(q\omega, \varphi^q)$.
  If $(\omega,\varphi)\in\cF$
and $\Phi:\torus\to\R$ is a lift of $\varphi$, then
$(\omega,\Phi):\kreis\times\R\to\kreis\times\R,\ (\theta,x)\mapsto
(\theta+\omega,x+\Phi(\theta,x))$ is a lift of $(\omega,\varphi)$. Thus, the
rotation numbers $\rho(\omega,\Phi,\theta)$ (or simply $\rho(\omega,\Phi)$ in
the quasiperiodic case) are defined by (\ref{e.fibred_rotnum}). Note that if $\omega=p/q$ is
rational, with $p,q$ coprime integers, then
$\rho(\omega,\Phi,\theta)$ is equal to $1/q$ times the rotation number of $F^q_\theta$, where $F=(\omega,\Phi)$.

Further, we will need two basic statements concerning the behaviour of rotation
numbers. Since both are folklore and follow from standard arguments, the proofs
are omitted. The first statement concerns the continuity of the boundary of
Arnold tongues.

\begin{lem} \label{l.arnold-boundaries}
  Suppose $\mathfrak{A}$ is a compact topological space and
  $(g_{\alpha,\tau})_{\alpha\in\mathfrak{A},\tau\in\R}$ is a two-parameter family of
  circle homeomorphisms (that is, $(\alpha,\tau,x)\mapsto g_{\alpha,\tau}(x)$ is
  continuous) with lifts $G_{\alpha,\tau}$ such that $(\alpha,\tau,x)\mapsto
  G_{\alpha,\tau}(x)$ is continuous and strictly monotone in $\tau$ for each
  $\alpha,x$ ({\em twist condition}).

  Then for any $\rho\in\R$ there exist continuous functions
  $\tau_\rho^-,\tau_\rho^+:\mathfrak A\to\R$ such that
  \[
    \{(\alpha,\tau)\in\mathfrak A\times\R \mid \rho(G_{\alpha,\tau})=\rho\} \ =
    \ \{(\alpha,\tau)\in\mathfrak A\times\R\mid \tau^-_\rho(\alpha)\leq \tau\leq
    \tau^+_\rho(\alpha)\} \ .
  \]
 \end{lem}
 The analogue statement even holds in the qpf case, see
 \cite{bjerkloev/jaeger:2009}.

 The second lemma concerns the uniform convergence of the fibred rotation number in
 parameter families.
\begin{lem}\label{l.rotnum-convergence}
  Suppose that $\cT$ is a compact topological space, $\omega\in\kreis$ is irrational,\\
  $\left(f_\tau= (\omega,\varphi_\tau) \right)_{\tau\in\cT}$ $ \subseteq\cF$ is a continuous family of
  qpf circle homeomorphisms with lifts $(F_\tau)_{\tau\in\cT}$,
  $\left(f_{n,\tau}=(\omega_n,\varphi_{n,\tau})\right)_{\tau\in\cT,n\in\N}$ is a
  sequence of such parameter families with lifts $(F_{n,\tau})_{\tau\in\cT,n\in\N}$ and $\nLim
  F_{n,\tau}=F_\tau$ for all $\tau\in\cT$. Then for all $\eps>0$ there exist
  $m_0,n_0\in\N$ such that for $m\geq m_0$ and $n\geq n_0$ and all
  $\tau\in\cT,\ \theta\in\kreis$ we have
  \[
       \left| \frac{1}{m}\left(F^m_{n,\tau,\theta}(x)-x\right) - \rho(F_\tau)\right| \ < \ \eps \ .
  \]
  In particular, $\rho(F_{n,\tau},\theta)\to \rho(F_\tau)$ uniformly on
    $\kreis\times\cT$ as $n\to\infty$.
\end{lem}

\section{Strategy of the proof}\label{Strategy}

Our key result, which entails the others rather easily (as explained in
Section~\ref{FurtherProofs}), is Theorem~\ref{t.main-homeomorphisms}.  For its
proof, we only have to show that mode-locked qpf circle homeomorphisms
$(\omega,\varphi)$ are dense in $\cF$, since the openness is obvious from
Definition~\ref{d.mode-locking}. In order to do so, we fix an initial map
$f=(\omega,\varphi)\in\cF$ and perturb it in several steps to obtain the
existence of an annulus that is mapped strictly inside itself by some
iterate. This will be done in four main steps.
\begin{enumerate}
\item[{\bf Step 1}] First, we replace the rotation number $\omega$ in the base
  by some rational numer $p/q$ to obtain a periodically forced circle
  homeomorphism.  By doing so, we may lose the uniqueness of the rotation
  number. However, by continuity ( see Lemma \ref{l.rotnum-convergence}), the rotation numbers $\rho(p/q,\Phi,\theta)$
  are still close to the original $\rho(\omega,\Phi)$. We then show that by
  means of a small modification of $\varphi$ we can ensure that the rotation
  numbers on the fibres are all the same and take a rational value. This is done
  by an indirect argument: we assume that no system in a neighbourhood of
  $(\omega,\varphi)$ is mode-locked (so that mode-locking is not dense) and build
  the required perturbation only under this assumption (since otherwise there is
  nothing to prove anyway).

\item[{\bf Step 2}] Once this is achieved, it follows from the standard theory
  of orientation-preserving circle homeomorphisms that all the fibre maps
  $f_\theta^q$ have a fixed point. We then proceed to show that $\varphi$
  can be perturbed in such a way that all the fibre maps $f^q_\theta$ are
  mode-locked.

  However, this does not yet imply that $f$ itself is mode-locked in the sense
  of Definition \ref{d.mode-locking} -- a forward invariant or periodic annulus
  $\cA$ does not yet have to exist for this perturbation. This happens, for
  example, if $\omega=p/q$ and the complement of set
  $B=\{(\theta,x)\in\torus\mid \Phi^q(\theta,x)=0\}$ is inessential.\footnote{An
    open set $U\subseteq \T^2$ is called {\em essential} if it contains some
    homotopically non-trivial closed curve $\gamma\subseteq U$, and it is called
    {\em inessential} otherwise. It is called doubly essential, if it contains
    two homotopically non-trivial closed curves of different homotopy types. An
    arbitrary set $A$ is called essential (doubly essential), if any open
    neighbourhood is essential (doubly essential).} In this case, the boundary
  of any homotopically non-trivial annulus $\cA\ssq \torus$ has to intersect
  $B$, and can therefore not be mapped into its own interior by any iterate of
  $f$ (since points of $B$ are fixed under $f^q$.

  What we obtain so far, though, is the fact that the sets
  $B^-=\{(\theta,x)\in\torus\mid \Phi^q(\theta,x)<0\}$ and
  $B^+=\{(\theta,x)\in\torus\mid\Phi^q(\theta,x)>0\}$ both project to all of
  $\kreis$.
\item[{\bf Step 3}] We then perturb $\varphi$ to ensure that both $B^-$ and
  $B^+$ consist of only finitely many connected components with positive
  distance to each other, and have some further good properties.
\item[{\bf Step 4}] The final step consists in showing that, in the situation obtained in Step
  3, there exists a small neighbourhood $B_\eps(p/q)$ such that for every
  $\omega'\in B_\eps(p/q)\smin\{p/q\}$ the system $(\omega',\varphi)$ is
  mode-locked in the sense of Definition~\ref{d.mode-locking}.
\end{enumerate}
While the arguments in Steps 1-3 are rather standard, it is the last step which
contains the main new geometric idea of our proof and depends crucially on the
topological characterisation of mode-locking in
Definition~\ref{d.mode-locking}. It is also this part in which the difference to
the linear-projective case of quasiperiodic \sltr-cocycles becomes most
apparent. What happens is that we have to deal with the fact that the desired
annulus $\cA$ may be of any homotopy type, whereas in the linear case the
corresponding forward invariant annuli (which correspond to invariant cone
fields for the cocycle) can only have homotopy types of the form $(1,k)$ with
$k\in\Z$ (that is, they wrap around the torus exactly once in the horizontal
direction).

\section{Proof of Theorem~\ref{t.main-homeomorphisms}}
\label{MainProofs}

As mentioned, we fix $f=(\omega,\varphi)\in\cF$ with lift $F=(\omega,\Phi)$ and
show that there exists an arbitrarily small perturbation of $f$ in $\cF$ that
exhibits mode-locking, following the steps outlined in the previous section. The
following proposition corresponds to Step 1.
\medskip

\begin{prop}\label{Step 1}
  Suppose that  $(\omega,\varphi)\in\cF$ and
  there exists a neighbourhood $B_\eta(\omega,\varphi)$ in $\cF$ such that no
  $f\in B_\eta(\omega,\varphi)$ is mode-locked. Then there exists
  $(\tilde\omega,\tilde\varphi)\in\cF$ such that
  \begin{itemize}
  \item $(\tilde\omega,\tilde\varphi)$ is an arbitrarily small perturbation of
    $(\omega,\varphi)$;
  \item $\tilde \omega=p/q$, where $p,q$ are coprime integers;
  \item there exists $m_0\in\N$ such that
    $\rho(\tilde\omega,\tilde\Phi,\theta)=m_0/q$ for all $\theta\in\kreis$.
  \end{itemize}
\end{prop}
\proof Fix $0<\eps<1$. We may assume that $\omega$ is irrational, so that
$(\omega,\Phi)$ has a unique rotation number. (Otherwise, we just replace
$\omega$ by some nearby irrational number). Moreover, since systems in a neighbourhood
of $(\omega,\varphi)$ are not mode-locked, there exists $\delta>0$ such that
\begin{equation} \label{e.S1_1}
  \rho(\omega,\Phi-\eps) + \delta \ < \ \rho(\omega,\Phi) \ <
  \ \rho(\omega,\Phi+\eps) - \delta\ .
\end{equation}
We now replace $\omega$ by $\tilde\omega=p/q$ with $1/q<\delta/2$. Then Lemma
\ref{l.rotnum-convergence} implies that if $\tilde\omega$ is chosen sufficiently
close to $\omega$, then
\begin{equation} \label{e.S1_2}
  \left|\rho(\tilde\omega,\Phi+\tau,\theta)-\rho(\omega,\Phi+\tau)\right|
  \ < \ \delta/2 \
\end{equation}
for all $\tau\in[-1,1]$. Fix $m_0\in\N$ such that
$|\rho(\omega,\Phi)-m_0/q|<\delta/2$. By Lemma~\ref{l.arnold-boundaries}, there
exists a continuous $1/q$-periodic function $\tau^+:\kreis\to\R$ such that
$\rho(\tilde\omega,\Phi+\tau^+,\theta)=m_0/q$ for all $\theta\in\kreis$. In order
to see this, apply Lemma \ref{l.arnold-boundaries} with $\mathfrak A=\kreis,\
\alpha=\theta$, $g_{\tau,\theta} =R_\tau\circ f_{\theta+(q-1)\tilde\omega}\circ
\ldots \circ R_\tau\circ f_\theta$ and $\rho=m_0$ and note that since the maps
$g_{\tau,\theta}$ and $g_{\tau,\theta+ip/q}$ are conjugate for all $i=1\ld q-1$,
the resulting functions $\tau^\pm=\tau^\pm_{m_0}$ are $1/q$-periodic. By
definition, we also have
$\rho(G_{\tau,\theta})=q\rho(\tilde\omega,\Phi+\tau,\theta)$, where $G_{\tau,\theta}=R_\tau\circ F_{\theta+(q-1)\tilde\omega}\circ \ldots\circ R_\tau\circ F_\theta$, so
$\rho(\tilde\omega,\Phi+\tau^+,\theta)=m_0/q$.

We claim that $|\tau^+|<\eps$. This holds since
\[
\rho(G_{\eps,\theta}) \ =
\ q\rho(\tilde\omega,\Phi+\eps,\theta)\ \stackrel{\eqref{e.S1_2}}{>}
\ q\rho(\omega,\Phi+\eps) - q\delta/2 \ \stackrel{\eqref{e.S1_1}}{>}
\ q\rho(\omega,\Phi)+q\delta/2 \ > \ m_0 \ ,
\]
and likewise $\rho(G_{-\eps,\theta})<m_0$ for all $\theta\in\kreis$.

In conclusion, since $\eps>0$ was arbitrary, we obtain the required perturbation
$(\tilde\omega,\tilde\varphi)=(p/q,\varphi+\tau^+)$.  \qed\medskip

As in Section \ref{Strategy}, if $\omega=p/q$ is rational with $p, q$ coprime and $(\omega,\varphi)$ satisfies the assertions of Proposition \ref{Step 1},
we let
\begin{eqnarray*}
B \ \ (\varphi) & = & \{(\theta,x)\in\T^2\ |\ \Phi^q(\theta,x)=0\} \ ,\\
B^-(\varphi)&=&\{(\theta,x)\in
\T^2\ |\ \Phi^q(\theta,x)<0 \},\\
B^+(\varphi)&=&\{(\theta,x)\in \T^2\ |\
\Phi^q(\theta,x)>0 \} \ ,
\end{eqnarray*}
where $\Phi^q$ is the lift of $\varphi^q$ such that $\{x\in\T^1\ |\ \Phi^q(\theta,x)=0\}\neq \emptyset$ for all $\theta\in\T^1$.
Then the following proposition takes care of Step 2.\medskip

\begin{prop}\label{Step 2}
  Suppose $(\omega,\varphi)\in \cF$ and there exists a neighborhood
  $B_\eta(\omega,\varphi)$ in $\cF$ such that no $f\in B_\eta(\omega,\varphi)$
  is mode-locked. Then there exists $\tilde f=(\tilde \omega,\tilde \varphi)\in
  \cF$ such that
\begin{itemize}
\item $(\tilde \omega, \tilde \varphi)$ satisfies the assertions of Proposition \ref{Step 1};
\item $\pi_1(B^+(\tilde\varphi))=\pi_1(B^-(\tilde\varphi))=\T^1$.
\end{itemize}
\end{prop}
Note that the last property means that all the fibre maps $\tilde f_\theta^q$ of
$\tilde f^q$ are mode-locked (as maps of the circle).\medskip

For the proof, we need two auxiliary lemmas that address the following technical
issue: If we perturb the map $f$ to obtain certain qualitative properties of the
fibre maps $f_\theta^q$, we need to take into account that these fibre maps are not
independent. More precisely, if $\theta, \theta'$ belong to the same periodic
orbit in the base, both $f_\theta^q$ and $f_{\theta'}^q$ result from the
composition of the same family of fibre maps $f_{\theta+\frac{i}{q}}$, and only
differ in the order of composition. Most importantly, they are conjugate to each
other. This is the content of the following

\begin{Lemma}\label{l.step2-conjugate}
  Suppose $\omega=p/q$ is rational, where $p, q$ are coprime integers. Then for
  any $\theta\in \T^1, j=1,2,\ldots, q-1$, the fibre maps $f_\theta^q$ and
  $f_{\theta+j/q}^q$ are conjugate.
\end{Lemma}
\proof Since $p,q$ are coprime, there exist unique $s,k\in\N$, such that
$sp-kq=1$, which means that for any $j=1,2,\ldots, q-1$, we have
$\frac{p}{q}\cdot js=jk+\frac{j}{q}$. Further
$f_\theta^{js+q}=f_{\theta+\frac{p}{q}\cdot js}^q\circ f_\theta^{js}=f_\theta
^{js}\circ f_\theta^q.$ Together with the fact that $f_{\theta+n}=f_\theta$ for
all $n\in\N$, this implies that
\[f_{\theta+\frac{j}{q}}^q=f_{\theta+\frac{p}{q}\cdot js}^q=f_\theta^{js}\circ
f_\theta^q\circ (f_\theta^{js})^{-1}.\]\qed\medskip

The next lemma allows to restrict to an interval $[\theta,\theta+1/q)$ when
perturbing $(p/q,\varphi)$.

\begin{Lemma}\label{l.step2-restr-interval}
  Let $\omega=p/q$ be rational with $p,q\in\N$ coprime, and $I\subseteq \T^1$ be
  a half open interval with $|I|\leq 1/q$. Suppose $\psi\in \Xi$ and
  $\psi(\theta,\cdot)=\varphi^q(\theta,\cdot)$ for $\theta\in\partial I$. Then
  there exists $(\omega, \tilde \varphi)\in \cF$ such that for any $\theta\in
  I$, we have $\tilde \varphi^q(\theta,\cdot)=\psi(\theta,\cdot)$.

  Moreover, $(\omega,\tilde \varphi)$ can be chosen as an arbitrary small
  perturbation of $(\omega,\varphi)$ provided that $\psi(\theta,\cdot)$ is
  sufficiently close to $\varphi^q(\theta,\cdot)$ for all $ \theta\in I$.
\end{Lemma}
\proof Let $f=(\omega,\varphi)$ and $g=(\omega,\psi)$ and define
\[\tilde \psi(\theta,x)=\left\{\begin{array}{ll}
\psi(\theta,x)\ & \textrm{for}\ \theta\in I,\ x\in\T^1\\
\varphi^q(\theta,x)\ & \textrm{for}\ \theta\in \T^1\setminus I,\ x\in\T^1.
\end{array}\right.\]
Then $\tilde \psi\in\Xi$ and hence $\tilde g=(\omega,\tilde \psi)\in \cF$. For
$(\theta,x)\in\T^2$, we let $\tilde\varphi(\theta,x)=\varphi_{\tilde
  f}(\theta,x)$ with $\tilde f_{\theta}(x)=f_{\theta}^{-(q-1)}\circ \tilde
g_\theta(x)$.
Then $\tilde f=(\omega, \tilde \varphi)\in\cF$ and we obtain
\begin{equation}\label{eq.step2-perturb-lemma}
\tilde f_\theta(x)=\left\{\begin{array}{ll}
f_\theta^{-(q-1)}\circ g_\theta(x)\ &\textrm{for}\ \theta\in I, x\in\T^1\\
f_\theta(x)\ & \textrm{for}\ \theta\in \T^1\setminus I, x\in\T^1.
\end{array}\right.
\end{equation}
Due to the fact that for any $\theta\in I, l=1,2,\ldots, q-1$, we have
$\theta+l\frac{p}{q}$ (mod 1) $\notin I$, we obtain
 \[\tilde f_\theta^q=\tilde f_{\theta+p/q}^{q-1}\circ
 \tilde f_\theta=f_{\theta+p/q}^{q-1}\circ f_\theta^{-(q-1)}\circ
 g_\theta= g_\theta,\ \forall \ \theta\in I\] and hence
 $\tilde\varphi^q(\theta,\cdot)=\psi(\theta,\cdot)$ for any $\theta\in I$ as required.

 Moreover, by (\ref{eq.step2-perturb-lemma}), for $\theta\in I$, we have $\tilde
 f_\theta-f_\theta=f_\theta^{-(q-1)}\circ g_\theta -f_\theta^{-(q-1)}\circ
 f_\theta^q$, which implies that $\tilde f$ is arbitrarily close to $f$ if
 $\psi(\theta,\cdot)$ is sufficiently close to $\varphi^q(\theta,\cdot)$ for all
 $\theta\in I$. \qed\medskip

 \proof[\bf{Proof of Proposition \ref{Step 2}}] Fix $\varepsilon>0$.  Without loss of
 generality, we may assume that $(\omega,\varphi)$ satisfies the assertions of
 Proposition \ref{Step 1}.  For $\theta\in\T^1$, we let
 \begin{eqnarray*}
   A(\theta,\Phi)& =& \{x\in\T^1\ |\ \Phi^q(\theta,x)=0\},\\
   A^+(\theta,\Phi) & = & \{x\in\T^1\ |\ \Phi^q(\theta,x)>0\},\\
 A^-(\theta,\Phi)&=&\{x\in\T^1\ | \ \Phi^q(\theta,x)<0\} \ .
\end{eqnarray*}
It is apparent that $\pi_1(B^+(\varphi))=\pi_1(B^-(\varphi))=\T^1$ is equivalent to
$A^+(\theta,\Phi)\neq \emptyset\neq A^-(\theta,\Phi)$ for all
$\theta\in\T^1$. Note that this also implies $A(\theta,\Phi)\neq \emptyset$ for
all $\theta\in\T^1$ by continuity.  We proceed in two steps and start by showing
that there exists $\hat \varphi\in B_{\eps/2}(\varphi)\cap \Xi$ such that
$A(\theta,\hat \Phi)$ is neither empty nor equal to $\T^1$ for all
$\theta\in\T^1$, which implies that $A^+(\theta,\hat \Phi)\cup A^-(\theta,\hat
\Phi)\neq \emptyset$ for all $\theta\in\T^1$ and hence
$\pi_1(B^-(\hat\varphi))\cup\pi_1(B^+(\hat\varphi))=\kreis$. Secondly, we show that $\hat \varphi$ can be
modified further and there exists $\tilde \varphi\in
B_{\eps/2}(\hat\varphi)\cap\Xi$ such that $\pi_1(B^-(\tilde\varphi))=\pi_1(B^+(\tilde\varphi))=\T^1$.\medskip

We may assume that there is some open interval of $\theta\in\T^1$ for which
$A(\theta,\Phi)\neq \T^1$, that is, $f_\theta^q\neq \mathrm{Id}_{\kreis}$.
Otherwise, using Lemma \ref{l.step2-restr-interval}, we can perturb
$(\omega,\varphi)$ in a small interval on the base such that $A(\theta,\Phi)\neq \T^1$ for
$\theta$ in the interior of this interval and the perturbed map still satisfies
the assertions of Proposition \ref{Step 1}. Without loss of generality, we assume that
$\theta=0$ lies in this interval and hence $A(\theta,\Phi)\neq \T^1$ for
$\theta=0$. Then by Lemma \ref{l.step2-conjugate}, $A(\theta,\Phi)\neq \T^1$ for
$\theta=1/q$, i.e. $f_{1/q}^q\neq \mathrm{Id}_{\kreis}$. Let
 \[\Theta\ = \  \left\{\theta\in[0,1/q]\mid  f_\theta^q=\mathrm{Id}_{\kreis} \right\}\ .\]
 Then, since $\Theta$ is compact and does not contain $0$ and $1/q$, we have
 that for any $\delta>0$ there exist finite number of disjoint closed intervals
 $J_1,\ldots, J_M$ and disjoint open intervals $U_1,\ldots, U_M$, such that
 $J_i\subseteq U_i\subseteq (0,1/q), (i=1,\ldots, M), \Theta\subseteq\cup_{i=1}^M
 J_i$, and for all $\theta\in \cup_{i=1}^M \textrm{cl}(U_i)$ and $x\in \T^1$, we
 have $|\varphi^q(\theta,x)|<\delta$. Then for each $U_i, (i=1,\ldots, M)$,
 there exists a continuous function $\eta_i: \textrm{cl}(U_i)\times
 \T^1\rightarrow \T^1$ which is differentiable with respect to $x$ and satisfies
\begin{itemize}
\item $|\eta_i(\theta,x)|<\delta, \forall\ \theta\in \textrm{cl} (U_i),
  x\in\T^1$, and $\partial _x \eta_i(\theta,x)>-1, \forall\ \theta\in U_i,
  x\in\T^1$;
\item for any $\theta\in U_i$, there exist $x_1,x_2\in \T^1$ such that
  $\eta_i(\theta,x_1)=0$ and $\eta_i(\theta,x_2)\neq 0$;
\item $\eta_i(\theta,x)=\varphi^q(\theta,x)$ for $\theta\in \partial U_i,
  x\in\T^1$.
\end{itemize}
Let
\[
\psi_1(\theta,x)\ = \ \left\{\begin{array}{ll} \eta_i(\theta,x) \quad \textrm{for}\ \theta\in U_i\\ \ \\
    \varphi^q(\theta,x) \quad \textrm{for}\ \theta\in \T^1\setminus \cup_{i=1}^M
    U_i\end{array} \ .\right.
\] Then $\psi_1\in\Xi$, and for all $(\theta,x)\in\T^2$, we have
$|\psi_1(\theta,x)-\varphi^q(\theta,x)|<2\delta$.  Applying Lemma
\ref{l.step2-restr-interval} with $I=[0,1/q), \psi=\psi_1$, there exists $\hat
\varphi\in \Xi$, such that $\hat \varphi^q(\theta,x)=\psi_1(\theta,x)$ for all
$\theta\in[0,1/q),\ x\in\T^1$. Due to the choice of $\eta_i$, this implies that
$A(\theta,\hat \Phi)$ is neither empty nor $\T^1$ for $\theta\in [0,1/q)$. Then
by Lemma \ref{l.step2-conjugate}, this holds for all $\theta\in\T^1$.  Moreover,
by Lemma \ref{l.step2-restr-interval}, if we choose $\delta$ small enough we can
ensure that $\hat f=(\omega,\hat \varphi)\in B_{\varepsilon/2}(\omega,\varphi)$.
This completes the first part of the proof. \medskip

Now, let $J, J'$ be two closed intervals of $\T^1$ with $|J|, |J'|<1/q$ and
\begin{equation}\label{c.step2-1}
\bigcup_{i=0}^{q-1}(J+i\omega)\cup(J'+i\omega)=\T^1.
\end{equation}
Choose two left closed right open intervals $\tilde J, \tilde {J}'\subseteq
\T^1$ with $|\tilde J|= |\tilde J'|=1/q$, such that $J\subseteq
\textrm{int}(\tilde J), J'\subseteq \textrm{int}(\tilde J')$. For
$\theta\in\T^1$, denote
\[\gamma^+_{\hat\Phi}(\theta)=\sup_{x\in\T^1}\hat \Phi^q(\theta,x),\ \
\textrm{and}\ \ \gamma_{\hat \Phi}^-(\theta)=-\inf_{x\in\T^1}\hat
\Phi^q(\theta,x).\] Then there exists $C>0$ such that for any $\theta\in\T^1$,
we have
\[\gamma_{\hat \Phi}^\pm(\theta)\in[0,C]
\quad \textrm{and}\quad \gamma_{\hat\Phi}^+(\theta)+\gamma_{\hat\Phi}^-(\theta)>0,\]
where the last inequality follows from the fact that
$\pi_1(B^-(\hat\varphi))\cup\pi_1(B^+(\hat\varphi))=\T^1$.  Moreover, there exists a continuous function
$\lambda: \T^1\rightarrow [0,1]$ supported on $\tilde J$ and $\lambda_{\mid
  J}=1$. For $0<\zeta<1/2$, we let
\[\psi_2(\theta,x)\ = \ \hat \varphi^q(\theta,x)-\zeta\lambda(\theta)\gamma_{\hat \Phi}^+(\theta)
+\zeta\lambda(\theta)\gamma_{\hat \Phi}^-(\theta) \] and apply
 Lemma \ref{l.step2-restr-interval} with $\psi=\psi_2,\ I=\tilde J$ and
 $\varphi=\hat \varphi$. Then there exists $\bar\varphi \in\Xi$, such that
 $\bar\varphi^q(\theta,x)=\psi_2(\theta,x)$ for $\theta\in\tilde J, x\in\T^1$.
 Therefore, for $\theta\in J$, we obtain that
 \[\gamma_{\bar\Phi}^+(\theta)=(1-\zeta)\gamma_{\hat \Phi}^+(\theta)+\zeta\gamma_{\hat\Phi}^-(\theta)>0.\]
In the same way, we obtain $\gamma_{\bar \Phi}^-(\theta)>0$ for $\theta\in J$.
Thus, we have $\gamma_{\bar \Phi}^\pm(\theta)>0$ for all
$\theta\in\cup_{i=0}^{q-1} J+i\omega$ by Lemma
\ref{l.step2-conjugate}. Moreover, as $d_{\cC^0}(\psi_2, \hat \varphi^q)\leq
2C\zeta$, we can ensure that $d_{\cC^0}(\bar \varphi, \hat
\varphi)<\varepsilon/4$ by choosing $\zeta$ is small enough, see
Lemma~\ref{l.step2-restr-interval}.

As for $J', \tilde J'$, similarly as before, we let
\[\psi_3(\theta,x)=\bar\varphi^q(\theta,x)-\tilde\zeta\tilde \lambda(\theta)
\gamma_{\bar\Phi}^+(\theta)+\tilde\zeta\tilde
\lambda(\theta)\gamma_{\bar\Phi}^-(\theta) \ , \] where $0<\tilde\zeta<1/2, \tilde \lambda:
\T^1\rightarrow [0,1]$ is continuous and supported on $\tilde J'$ and $\tilde
\lambda_{|J'}=1$. Applying Lemma~\ref{l.step2-restr-interval} once more, with
$\psi=\psi_3, I=\tilde J', \varphi=\bar\varphi$, we obtain a
$\tilde\varphi\in\Xi$ such that $\tilde \varphi^q(\theta,x)=\psi_3(\theta,x)$
for $\theta\in\tilde J', x\in\T^1$. Then,  we have
\begin{equation}\label{e.step2-gamma-differ}
\gamma_{\tilde \Phi}^\pm(\theta)>0
\end{equation}
for all $\theta\in \bigcup_{i=0}^{q-1} J'+i\omega$. It remains to show that (\ref{e.step2-gamma-differ}) now
holds for all $\theta\in\kreis$. To that end, suppose $\theta\in \tilde
J'\setminus \bigcup_{i=0}^{q-1} J'+i\omega$. Then by the choice of $J, J', \tilde J'$, we have that
$\theta\in \cup_{i=0}^{q-1}(J+i\omega)$. Therefore,
\[\gamma_{\tilde \Phi}^+(\theta)\ = \ (1- \tilde\zeta\tilde\lambda(\theta))
\gamma_{\bar\Phi}^+(\theta)+\tilde\zeta\tilde \lambda(\theta)\gamma_{\bar
  \Phi}^-(\theta) \ \geq \ (1-\tilde \zeta)\gamma_{\bar\Phi}^+(\theta)\ > \ 0
\ .\] In the same way, we get $\gamma_{\tilde \Phi}^-(\theta)>0$ for $\theta\in
\tilde J'\setminus \bigcup_{i=0}^{q-1} J'+i\omega$. Hence, we have that (\ref{e.step2-gamma-differ}) holds
for all $\theta\in \tilde J'$. By Lemma \ref{l.step2-conjugate}, however, this
means that $\gamma_{\tilde \Phi}^\pm(\theta)>0$ for all $\theta\in \T^1$, and
thus $\pi_1(B^+(\tilde\varphi))=\pi_1(B^-(\tilde\varphi))=\T^1$. Moreover, as
$d_{\cC^0}(\psi_3,\bar\varphi^q)\leq 2C'\tilde \zeta$, where $C'=
\sup_{\theta\in\T^1}\{\gamma_{\bar \Phi}^+(\theta),
\gamma_{\bar\Phi}^-(\theta)\}<\infty$, Lemma~\ref{l.step2-restr-interval}
implies that if $\tilde \zeta$ is small enough, then we have $d_{\cC^0}(\tilde
\varphi, \bar \varphi)<\varepsilon/4$.  This completes the proof.  \qed\medskip

As mentioned before, even if $f=(\omega,\varphi)\in \cF$ satisfies the assertions of
Proposition \ref{Step 2} and all the fibre maps $f_\theta^q$ are mode-locked, this
does not imply the existence of some topological annulus $\cA$ that is mapped
inside itself by $f^q$. The reason is that possibly neither of $B^+(\varphi)$ or $B^-(\varphi)$ is
essential (meaning that $B(\varphi)$ is doubly essential). As explained in the preceeding
section, this is an obstruction to the existence of an annulus $\cA$ that is
mapped strictly inside itself by some iterate. Hence, we will perturb $\varphi$
to ensure that both $B^+$ and $B^-$ consist of finitely many connected
components and the distance between each two components of $B^+$ (resp. $B^-$)
is positive. This implies immediately that at least one of the two sets is
essential. The existence of such a perturbation is established in the following
proposition, which corresponds to Step 3. In order to state it, we need to introduce
some further notions.

Suppose that $\tilde B:=B(\tilde\varphi)=\{(\theta,x)\in \torus \mid \tilde\Phi^q(\theta,x)=0\}$
consists of a finite number of line segments. Call the endpoints of these line
segments {\em vertices} and further assume that each vertex is the endpoint of
exactly two line segments. Then the connected components of $\tilde B$ are
compact one-dimensional topological manifolds, and thus topological circles. We
say a vertex $v=(\theta,x)$ is {\em critical}, if there exists $\eps>0$ such
that $\pi_1(B_\eps(v)\cap \tilde B)$ is a half-sided interval with $\theta$ as
one endpoint.  We say $v$ is a {\em left critical vertex} if $\theta$ is the
left endpoint of that interval, and a {\em right critical vertex} if $\theta$ is
the right endpoint. Moreover, we use the same notions for vertices of the lift
$\cB$ of $\tilde B$ to $\R^2$. Similar to before, we let $\tilde
B^+:=B^+(\tilde\varphi)=\{(\theta,x)\in \torus \mid \tilde\Phi^q(\theta,x)>0\}$ and $\tilde
B^-:=B^-(\tilde\varphi)=\{(\theta,x)\in \torus \mid \tilde\Phi^q(\theta,x)<0\}$.\medskip

\begin{prop}\label{Step 3}
  Suppose that $(\omega,\varphi)\in \cF$ and there exists a neighborhood
  $B_\eta(\omega,\varphi)$ in $\cF$ such that no $f\in B_\eta(\omega,\varphi)$
  is mode-locked. Then there exists $\tilde f=(\tilde \omega,\tilde \varphi)\in
  \cF$ such that
\begin{itemize}
\item[(i)] $(\tilde \omega, \tilde \varphi)$ satisfies the assertions of Proposition
  \ref{Step 2};
\item[(ii)] $\tilde B$ consists of a finite number of line segments, and each
  vertex belongs to exactly two different line segments;
\item[(iii)] no line segment in $\tilde B$ is vertical;\footnote{Here we say a
    line segment on the torus is vertical if it lies on some fibre
    $\{\theta\}\times\T^1$}
\item[(iv)] two different connected components of $\tilde B^+$ have positive
  distance to each other, and the same is true for $\tilde B^-$;
\item[(v)] every fibre $\{\theta\}\times\kreis$ contains at most one critical
  vertex.
\end{itemize}
\end{prop}
Note that due to (ii), the connected components of $\tilde B$ are topological
circles. As a consequence, connected components of $\tilde B^-$ and $\tilde B^+$ are either
topological disks or annuli, or punctured topological disks or annuli.

\proof Fix $\eps>0$. Without loss of generality, we may assume that
$(\omega,\varphi)$ satisfies the assertions of Proposition \ref{Step 2}. In particular
$\omega=p/q$ is rational. We will not change the rotation number in the base and
thus set $\tilde \omega=\omega$.

In order to find the desired approximation $\tilde \varphi$, we work on the
level of lifts and construct an approximation $\tilde\Phi$ of $\Phi$, where
$\Phi:\torus\to\R$ is a lift of $\varphi$. To that end, we first fix a
triangularisation $\Delta=\{D_i\mid i=1\ld m\}$ of $\torus$ such that
\begin{equation} \label{e.triangle_size}
\diam(D_i) \ < \ \viertel \min_{k=1}^{q-1} d(0,k\omega) \ .
\end{equation}
Further, we require that all the vertices of the triangles in $\Delta$ and their
first $q-1$ iterates lie in pairwise different fibres. Note that the latter only
depends on the projections of the vertices to the $\theta$-coordinate and on the
rotation number $\omega$, but not on the fibre dynamics.  Denote by $v^j_i$ with
$j=1,2,3$ the three vertices of $D_i$. Obviously, we may have
$v^j_i=v^{j'}_{i'}$ for $(i,j)\neq(i',j')$, but every vertex belongs only to a
finite number of triangles and the diameter of the union of all triangles
containing a given vertex is smaller than $\min_{k=1}^{q-1} d(0,k\omega)/2$.

Now, we first approximate $\Phi$ by a continuous map $\tilde\Phi$ which is
piecewise affine on every triangle $D_i\in\Delta$. This means that, $\tilde\Phi$
is uniquely determined by the values $a^j_i=\tilde\Phi(v^j_i)$ (where
$a^j_i=a^{j'}_{i'}$ whenever $v^j_i=v^{j'}_{i'}$). Moreover, if $\sup_{i=1}^m
\diam(D_i)$ is sufficiently small and $d(a^j_i,\Phi(v^j_i))<\delta$ for all
$i=1\ld m,\ j=1,2,3$ and some sufficiently small $\delta$, then by uniform
continuity $\tilde\Phi$ will be $\eps$-close to $\Phi$. In all of the following,
we assume that this is the case, which is to say that we choose the
triangularisation sufficiently fine and pick all of the $a^j_i$ from
$B_\delta(\Phi(v^j_i))$.

Let $\tilde\varphi=\pi\circ\tilde\Phi$, $\tilde f=(\omega,\tilde\varphi)$ and let
further
\[
\cP_{\tilde\varphi} \ = \ \bigvee_{n=0}^{q-1} \tilde f^{-n}(\Delta) \
\]
be the common refinement of the partitions $\tilde f^{-n}(\Delta)$ (slightly
abusing standard terminology, we allow intersections between the partition
elements on their boundaries). Then, since $\tilde \Phi$ is piecewise affine,
$\cP_{\tilde\varphi}$ consists of polygons and $\#\cP_{\tilde\varphi}$ is a finite
number that is bounded by some number only depending on $m$ and $q$. Denote by $\vertset(\cP_{\tilde\varphi})$
the set of vertices of polygons in $\cP_{\tilde\varphi}$. We now show how the
$a^j_i$ can be adjusted to ensure (ii) and (iii), whereas (iv) follows by
construction and (v) can be achieved by some extra modification.
\begin{itemize}
\item[(ii)] We claim that the $a^j_i$ can be chosen such that
 \begin{equation}\label{e.nonzero_vertices}
   \tilde\Phi^q(v)\ \neq \ 0 \quad \textrm{for all } v\in\vertset(\cP_{\tilde\varphi}) \ .
 \end{equation}
 If this holds, then for all $P\in\cP_{\tilde\varphi}$ the graph of
 $\tilde\Phi^q_{|P}$, which is the restriction of a hyperplane to the region above
 $P$, cannot be contained in $\R^2\times \{0\}$. Moreover, the intersection of
 the hyperplane with $\R^2\times\{0\}$, which is either empty or a straight
 line, cannot contain any face or vertex of the polygon. Hence, the intersection
 of the graph of $\tilde\Phi^q_{|P}$ with $\R^2\times\{0\}$, which projects to
 $\tilde B\cap P$, consists of at most finitely many line segments (note that $P$ may
 not be convex, so that multiple line segments are possible). Further, the
 relative interiors of these line segments are contained in the interior of $P$,
 and their endpoints are contained in the relative interior of two faces of
 $P$. Since each of these faces belong to exactly one other polygon, each
 endpoint of a line segments in $\tilde B\cap P$ is contained in exactly one further
 line segment in $\tilde B$ that lies in the adjacent polygon.

 Since the total number of vertices in $\vertset(\cP_{\tilde\varphi})$ is bounded
 by some $N\in\N$ depending only on $m$ and $q$, condition
 (\ref{e.nonzero_vertices}) can be ensured by a successive application of the
 following claim.
 \begin{claim}
   Let $M(\cP_{\tilde\varphi})= \{v\in\vertset(\cP_{\tilde\varphi})\mid
   \tilde\Phi^q(v)\neq 0\}$. If $\#M(\cP_{\tilde\varphi}) = \ell$ and
   $\vertset(\cP_{\tilde\varphi})\smin M(\cP_{\tilde\varphi})$ is non-empty, then an
   arbitrary small modification of one of the values $a^j_i$ allows to obtain a
   perturbation $\bar\varphi$ of the same form as $\tilde\varphi$ that satisfies
   $\#M(\cP_{\bar\varphi}) \geq \ell+1$.
 \end{claim}
 \proof[Proof of the claim.] We first note that if $\bar\varphi$ is obtained
 from $\tilde \varphi$ by a sufficiently small modification of one of the $a^j_i$,
 then for all vertices $v\in M(\cP_{\tilde\varphi})$ there will be a vertex $\bar
 v\in M(\cP_{\bar\varphi})$ arbitrarily close to $v$ that still satisfies
 $\bar\Phi^q(\bar v)\neq 0$. Thus, $\#M(\cP_{\bar\varphi})\geq
 \#M(\cP_{\tilde\varphi})$ and it suffices to show that if
 $\vertset(\cP_{\tilde\varphi})\smin M(\cP_{\tilde\varphi})$ is non-empty, we may
 choose $v^*\in \vertset(\cP_{\tilde\varphi})\smin M(\cP_{\tilde\varphi})$ and a
 perturbation $\bar\varphi$ obtained from $\tilde\varphi$ by an arbitrarily small
 modification of one of the $a^j_i$ such that $v^*$ is also contained in
 $M(\cP_{\bar\varphi})$ but $\bar\Phi^q(v)\neq 0$.

 In order to do so, fix $v^*=(\theta_*,x_*) \in \vertset(\cP_{\tilde\varphi})\smin
 M(\cP_{\tilde\varphi})$. Note that by definition of $\cP_{\tilde\varphi}$ we must
 have $v^*\in \tilde f^{-k}(\partial D_r)$ for some $D_r\in\Delta$ and $k\in\{0\ld
 q-1\}$. Let $v'=(\theta',x')=\tilde f^{q-1}(v^*) \in D_i$ for some $i\in\{1\ld m\}$
 and suppose that $\bar \varphi$ is obtained from $\tilde \varphi$ by modifying
 $a^j_i$ at some vertex $v^j_i$ of $D_i$. Then, due to (\ref{e.triangle_size})
 and $\omega=p/q$, we do not modify $\bar\varphi(\theta,x)$ for any $\theta$ in
 one of the intervals $B_\eta(\theta_*)+n\omega$ with $\eta=\min_{k=1}^{q-1}d(0, k\omega)/2,\
 n=0\ld q-2$. However, this
 means that none of the maps $\tilde f\ld \tilde f^{q-1}$ is modified in a
 neighbourhood of $v^*$, and by definition of $\cP_{\tilde\varphi}$ and
 $\cP_{\bar\varphi}$, respectively, this means that $v^*\in \vertset(\cP_{\bar\varphi})$.
 Moreover, we
 have that $\bar\Phi^q(v^*) = \bar F_{\theta'}(x')-x_*$, since
 $\bar f^{q-1}_{\theta_*}(x_*)=\bar f_{\theta_*+(q-2)\omega}\circ\cdots\circ \bar f_{\theta_*}(x_*) =
 \tilde f^{q-1}_{\theta_*}(x_*)$. Thus, $\bar\Phi^q(v^*)$ depends monotonically on $\bar
 F_{\theta'}(x')$. The latter depends monotonically on the values $a^j_i$, and
 for at least one $j\in\{1,2,3\}$ this dependence is strictly monotone (recall
 that $\bar F$ is affine on $D_i$, so that $\bar F_{|D_i}$ is determined by the choice of
 the $a^j_i$). Hence, it is possible to modify $a^j_i$ slightly to ensure that
 $\bar\Phi^q(v)\neq 0$. This proves the claim.
  \qed
 \end{itemize}
 \begin{itemize}
\item[(iii)] In order to exclude the existence of vertical line segments in $\tilde B$, it
  suffices to ensure that $\partial_x \tilde\Phi_{|\textrm{int}(P)}^q\neq 0$ for all
  $P\in\cP_{\tilde\varphi}$. Similar as in (ii), this is a consequence of the
  following
  \begin{claim}
    Let
    $V(\cP_{\tilde\varphi})=\{P\in\cP_{\tilde\varphi}\mid \partial_x\tilde\Phi^q_{|\inte(P)}
    \neq 0\}$. If $\#V(\cP_{\tilde\varphi})=\ell$ and $\cP_{\tilde\varphi}\smin
    V(\cP_{\tilde\varphi})$ is non-empty, then an arbitrary small modification of
    one of the values $a^j_i$ allows to obtain a perturbation $\bar\varphi$ of
    the same form as $\tilde\varphi$ that satisfies $\#V(\cP_{\bar\varphi})\geq
    \ell+1$.
  \end{claim}
  \proof[Proof of the claim.] As in the previous claim, if $\bar\varphi$ is
  sufficiently close to $\tilde\varphi$, then the polygons in
  $V(\cP_{\tilde\varphi})$ will persists (possibly in slightly perturbed form) and
  still belong to $V(\cP_{\bar\varphi})$. Hence, we choose some
  $P\in\cP_{\tilde\varphi}\smin V(\cP_{\tilde\varphi})$. Then $\tilde f^{q-1}(P)$ is
  contained in some $D_i\in\Delta$. As before, if we modify $a^j_i$ at one of
  the vertices $v^j_i$ of $D_i$, then we will not alter $\tilde f$ on any of the
  polygons $P,\tilde f(P)\ld \tilde f^{q-2}(P)$ (note that the horizontal extent of these
  polygons is always the same and smaller than that of $D_i$). Hence, $P$
  remains unchanged, and the slope of $\bar\Phi^q$ in the $x$-direction on
  $P$ depends strictly monotonically on at least one of the values $a^j_i$. This
  allows to perform the required modification to render $\partial_x
  \bar\Phi^q_{|\textrm{int}(P)}\neq 0$.\qed
\end{itemize}
\begin{itemize}
\item[(iv)] By construction, whenever we cross a line segment in $\tilde B$, the sign
  of $\tilde\Phi^q$ changes. Hence, two different connected components of $\tilde B^-$,
  respectively $\tilde B^+$, cannot lie directly next to each other, so there is a
  minimal distance between two connected components of the same set.
\item[(v)] Finally, if there are two critical vertices in $\vertset(\tilde B)$ on the
  same fibre $\{\theta\}\times\kreis$, then we can use Lemmas
  \ref{l.step2-conjugate} and \ref{l.step2-restr-interval} again to perform a
  perturbation on $B_\epsilon(\theta)\times\kreis$, where $\epsilon>0$ is
  sufficiently small and such that there are no critical vertices of $\tilde B$ in
  $(B_\epsilon(\theta)\setminus\{\theta\})\times\T^1$, that moves these critical
  vertices to different fibres. Note here that critical vertices are always
  mapped to critical vertices (which may not be true for arbitrary
  vertices).\qed
\end{itemize}\medskip

It remains to show that any $(\tilde\omega, \tilde \varphi)\in\cF$ satisfying
the assertions of Proposition \ref{Step 3} is either mode-locked, or we can find some
perturbation of it that is mode-locked.  This is the content of the following
proposition, which corresponds to Step 4 in Section~\ref{Strategy} and completes
the proof of Theorem \ref{t.main-homeomorphisms}.

Before we turn to the precise statement and proof, we briefly comment on the
strategy. We assume that a perturbation $\tilde f=(\tilde\omega,\tilde\varphi)$
of our original system satisfies the assertions of Proposition~\ref{Step 3}. Then, in
the simplest case, both sets $\tilde B^-$ and $\tilde B^+$ may be essential and
contain two homotopically non-trivial closed curves $\gamma^-$ and $\gamma^+$,
respectively. In this situation, $\gamma^-$ is disjoint from its image and
mapped in the clockwise direction by $\tilde f^q$, whereas $\gamma^+$ is mapped
in the counterclockwise direction. As a consequence, the annulus which is
bounded below by $\gamma^+$ and above by $\gamma^-$ is mapped strictly inside
itself, so that $\tilde f$ is mode-locked.  The main problem in the general case
is that one of the two sets $\tilde B^\pm$, say $\tilde B^-$, may not be
essential. In this case, the boundary of any homotopically non-trivial annulus
$\cA$ has to intersect $\tilde B$. As points in $\tilde B$ are fixed by $\tilde
f^q$, this means $\cA$ cannot be mapped strictly inside itself.

The idea is now to construct a homotopically non-trivial closed curve $\gamma^-$
which contains some vertical segments, and while $\gamma^-$ is allowed to
intersect $\tilde B$ and also $\tilde B^+$, this happens only inside these
segments. Hence, all the non-vertical parts of $\gamma^-$ will be contained in
$\tilde B^-$ and therefore be mapped in the clockwise direction. Then the
vertical segments will still intersect their images, but they can be rendered
disjoint from their images by slightly changing the rotation number $\tilde
\omega$ in the base. If the vertical and non-vertical parts of $\gamma^-$ are
then arranged in the right way, one can ensure that the whole curve $\gamma^-$
is disjoint from its image and mapped in the clockwise direction. This will then
also remain true for a small modification of $\gamma^-$ in which the vertical
segments are slightly tilted, so that locally $\gamma^-$ can be represented as
the graph of a function over the base $\kreis$. The construction of $\gamma^+$
works in an analoguous way, and in fact, both curves have to be constructed
simultaneously in order to ensure that the homotopy types coincide.  \medskip

\begin{prop}\label{Step 4}
  Let $\tilde f=(\tilde\omega,\tilde \varphi)\in\cF$ satisfy the assertions of
  Proposition~\ref{Step 3}.  Then either $\tilde f$ is mode-locked, or there exists
  $s>0$ such that for any other rotation number $\omega'\in
  (\tilde\omega,\tilde\omega+s)$ the map $(\omega', \tilde \varphi)\in\cF$ is
  mode-locked.
\end{prop}
Before we turn to the proof, we first need to introduce some further notation.
Let $\cB, \cB^\pm$ be the lifts of $\tilde B, \tilde B^\pm$ in $\R^2$.
We denote by $\operatorname{VP}$ the set of vertices of $\cB$, and by \CV\ and
\RCV\ the sets of critical vertices and right critical vertices, respectively.
For any
triple $(\hat\theta,\hat x,\hat y)\in\R^3$, we say $(\hat \theta, \hat y)$ is in
the shadow of $(\hat \theta, \hat x)$ (or $(\hat \theta, \hat x)$ shadows $(\hat
\theta, \hat y)$) if $\hat y<\hat x$, $(\hat \theta,\hat x)\in \cB^-$, $(\hat
\theta,\hat y)\in \cB^+$, and there is exactly one $\hat \xi\in(\hat y, \hat
x)\subseteq\R$ with $(\hat \theta,\hat \xi)\in\cB$.  Then, by the assertions
of Proposition \ref{Step 3}, if $(\{\hat\theta\}\times\R^1) \cap CV=\emptyset$ and
$(\hat \theta, \hat x)\in\cB^-$ there always exists some $(\hat \theta, \hat
y)\in\cB^+$ in the shadow of $(\hat \theta, \hat x)$. We define
\[\cS \ = \ \left\{(\hat \theta, \hat x, \hat y)\in \R^3\ |\ (\hat \theta, \hat y)\
 \textrm{is in the shadow of}\ (\hat \theta, \hat x)\right\}.\] For $(\hat
 \theta, \hat x, \hat y)\in \cS$ and $\hat\theta'\geq \hat\theta$, let
\begin{eqnarray*}\mathcal{CS}(\hat\theta,\hat x,\hat y,\hat \theta') \ =
\ \left\{(\Lambda^-,\Lambda^+)\ \left|\ \begin{array}{lll}\Lambda^\pm :
       [\hat\theta,\hat\theta'] \rightarrow \R\ \textrm{are
           continuous},\\ \ \\ \Lambda^-(\hat \theta)=\hat x,\ \Lambda^+(\hat
         \theta)=\hat y,\ \textrm{ and}\\\ \\ (\hat \vartheta,
         \Lambda^-(\hat\vartheta),\Lambda^+(\hat\vartheta))\in\cS, \ \forall
         \ \hat\vartheta \in[\hat\theta,\hat\theta')\end{array}\right.\right\} \ .
\end{eqnarray*}
Further, we let
\[\hat \theta^*(\hat \theta,\hat x,\hat y)\ =\
\sup\left\{\hat
\theta'\in\R\ |\ \exists\ (\Lambda^-,\Lambda^+)\in\mathcal{CS}(\hat\theta,\hat
x,\hat y,\hat\theta')\right\} \ .\] Obviously $\hat\theta^*(\hat \theta,\hat
x,\hat y)\geq \hat\theta$. Moreover, the following statement holds.\medskip
\begin{Lemma}\label{l.target-point}
Let $(\hat \theta,\hat x,\hat y)\in\cS$. Then $\hat\theta^*=\hat
\theta^*(\hat\theta,\hat x,\hat y)>\hat\theta$.  Furthermore,
\begin{itemize}
\item if $\hat\theta^*=\infty$, then there exist closed curves
  $\gamma^\pm\subseteq \tilde B^\pm$ and they are in the same homotopy
  class. Moreover, the lifts $\Gamma^\pm$ of $\gamma^\pm$ in $\R^2$ are both
  simple curves;
\item otherwise, the set $\mathcal{CS}(\hat\theta,\hat x,\hat y,\hat\theta^*)$
  is non-empty and there exists $\hat x^*\in\R$ with $(\hat\theta^*,\hat
  x^*)\in RCV$ such that for any
  $(\Lambda^-,\Lambda^+)\in\mathcal{CS}(\hat\theta,\hat x,\hat
  y,\hat\theta^*)$ either $ \Lambda^-(\hat\theta^*)=\hat x^*$ or
  $\Lambda^+(\hat\theta^*)=\hat x^*$.
\end{itemize}
\end{Lemma}
\proof Let $\theta=\pi(\hat\theta)$ and $\tilde B_{\theta}^-=(\{\theta\}\times
\T^1)\cap \tilde B^-=\cup_{i=1}^{N_{\theta}}J_{\theta}^i$ where $J_\theta^i$
are the connected components of $\tilde B^-_\theta$. As $\tilde B$ consists only
of a finite number of line segments and $\partial \tilde B_\theta^-\ssq \tilde
B$, we have $N_\theta<\infty$. If $\hat \theta^*=\infty$, then there exist
simple curves $(\Lambda^-,\Lambda^+)\in\mathcal{CS}(\hat\theta,\hat x,\hat
y,\hat\theta+N_\theta)$.  Hence, there exist $i,j\in\{1,\ldots, N_{\theta}\}$
such that $\pi\big(\Lambda^-(\hat\theta+i)\big)$ and
$\pi\big(\Lambda^-(\hat\theta+j)\big)$ belong to the same component of
$\tilde B_{\theta}^-$. Without loss of generality we assume $i=0$. Then we can find
some $k\in\Z$ such that $\hat x+k$ and $\Lambda^-(\hat\theta+j)$ belong to the
same component of $\cB^-$. Since $\cB^-$ is open and contains the whole vertical
line segment between these two points, it is possible to modify $\Lambda^-$ in a
neighbourhood of $\hat\theta+j$ in order to obtain $\Lambda^-(\hat\theta+j)=\hat
x+k$ without leaving $\cB^-$.  Moreover, by the definition of shadow points, we
know that at the same time $\Lambda^+(\hat\theta+j)$ and $\hat y+k$ belong to the
same component of $\cB^+\cap \{\hat\theta+j\}\times\R$, and we can modify
$\Lambda^+$ such that $\Lambda^+(\hat\theta+j)=\hat y+k$ without leaving $\cB^+$.
Then the projections $\gamma^\pm$ of $\Lambda^\pm$ in $\T^2$ are closed curves in
the same homotopy class with $\gamma^\pm\subseteq \tilde B^\pm$ as required. \medskip

Now we consider the case when $\hat\theta^*<\infty$ and let
$\theta^*=\pi(\hat\theta^*)$. Further, we define
$\CV_{\hat\theta^*}=(\{\hat\theta^*\}\times \R)\cap \CV$ and
$\RCV_{\hat\theta^*}=(\{\hat\theta^*\}\times \R)\cap \RCV$ and denote by
$\LS_{\hat\theta^*}$ the set of line segments $\sigma$ of $\cB$ with
$\hat\theta^*\in \pi_1(\sigma)$. Due to $\#\pi(\VP)<\infty$, there exists $\eta_0>0$
such that
\begin{equation}\label{no-other-vertex}
\Big(\big(B_{\eta_0}(\hat\theta^*)\setminus\{\hat\theta^*\}\big)\times \R^1 \Big)\cap \VP=\emptyset.
\end{equation}
 Then for any $\sigma\in \LS_{\hat\theta^*}$, each endpoint of $\sigma$ either
 belongs to $\{\hat\theta^*\}\times \R$ or is outside
 $B_{\eta_0}(\hat\theta^*)\times\R$.  Due to choice of $\hat y$ in the shadow of $\hat
 x$, this makes it easy to see that $\hat\theta^*>\hat \theta$.
\begin{figure}[ht!]
 \epsfig{file=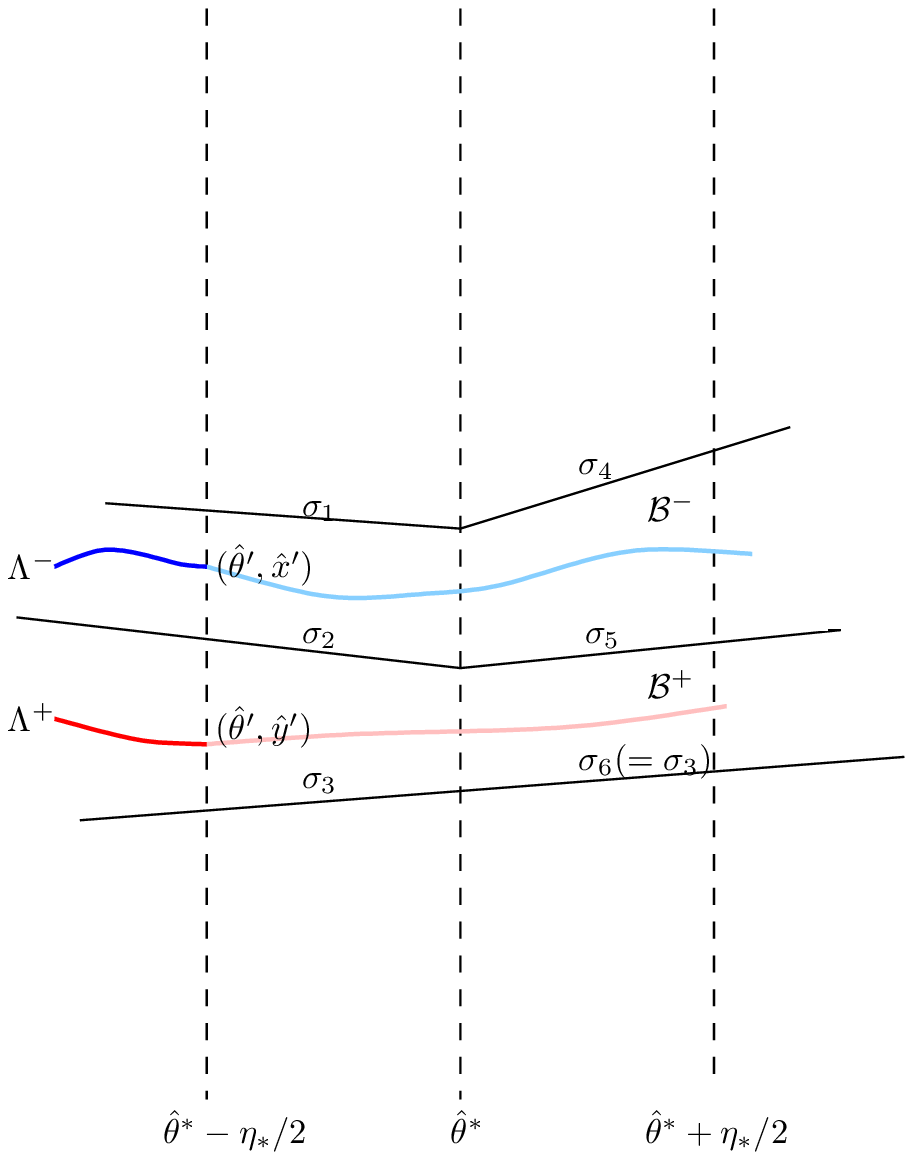, clip=,width=0.3\linewidth}\hspace{2eM}
 \epsfig{file=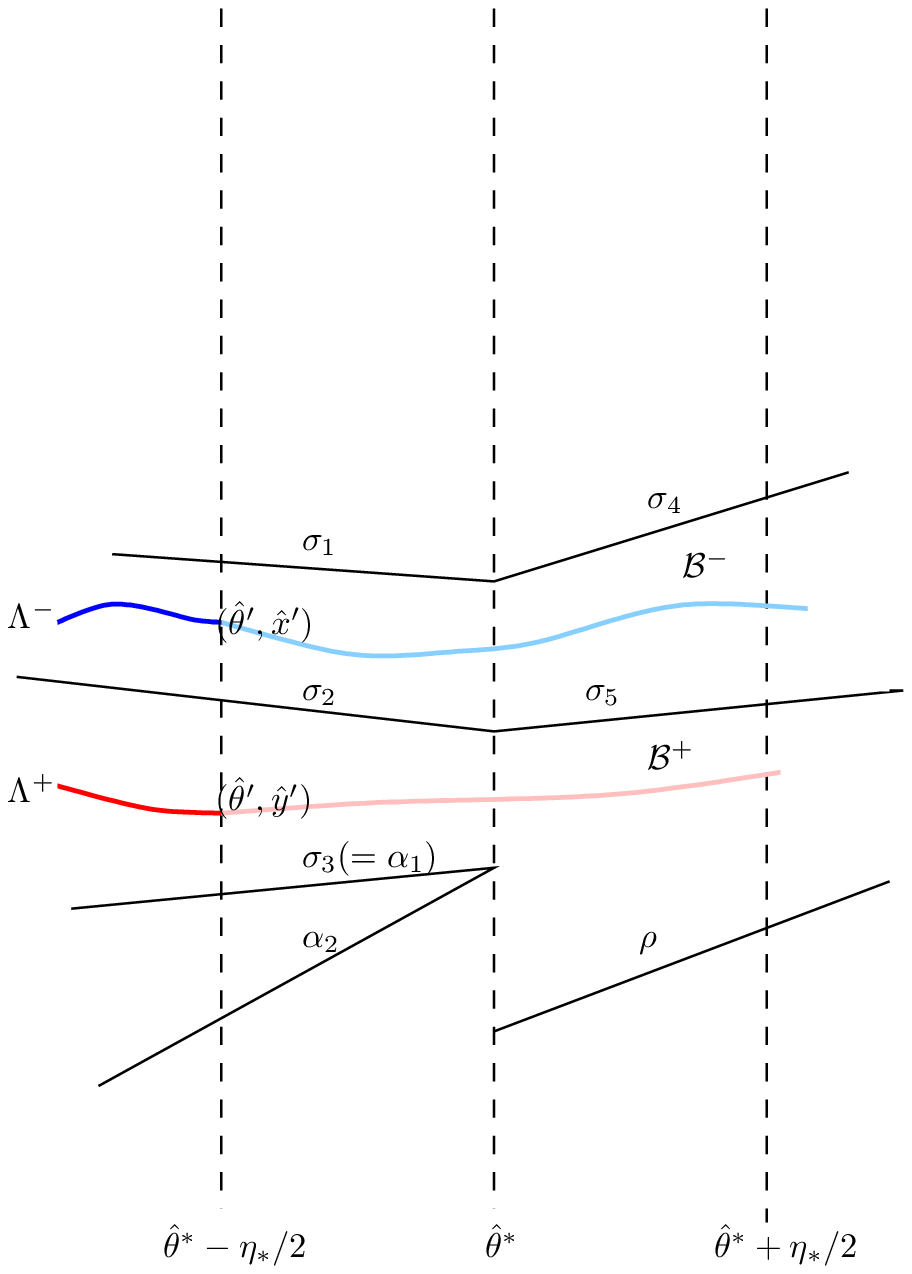,
   clip=,width=0.3\linewidth} \hspace{2eM}
 \epsfig{file=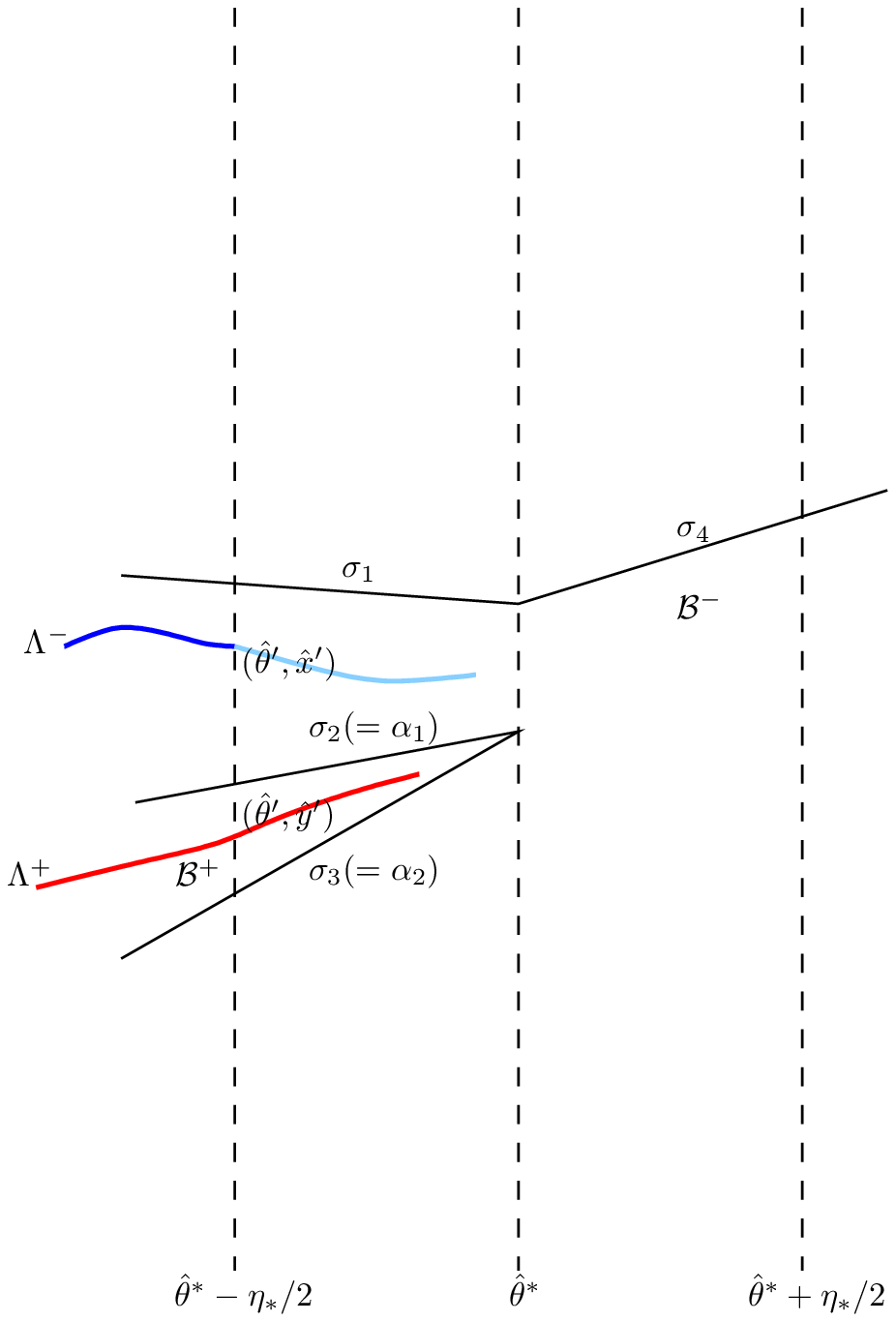, clip=,width=0.3\linewidth}

\vspace{2ex}
\hspace{5eM}  (a)\hspace{13eM}  (b)
\hspace{13eM}  (c)

\caption{Possible configurations for the line segments in
  $\LS_{\hat\theta^*}$. \label{Figure-extension} }
\end{figure}

Now, let $\eta_*=\min\{\hat\theta^*-\hat\theta, \eta_0\}>0$.  Then by the
definition of $\hat\theta^*$ there exist continuous curves $\Lambda^\pm:
[\hat\theta, \hat \theta^*-\eta_*/2]\rightarrow \R$ with $\Lambda^-(\hat
\theta)=\hat x,\ \Lambda^+(\hat \theta)=\hat y$ and
\[
(\hat\vartheta,\Lambda^-(\hat\vartheta),\Lambda^+(\hat\vartheta))\ \in \ \cS \quad \forall
\ \hat\vartheta \in[\hat\theta,\hat\theta^*-\eta_*/2] \ .
\]
Let $(\hat {\theta}',\hat { x}',\hat{y}'):=(\hat\theta^*-\eta_*/2,
\Lambda^-(\hat\theta^*-\eta_*/2), \Lambda^+(\hat\theta^*-\eta_*/2))\in\cS$ and
$\theta'=\pi(\hat{\theta}')$. Then, due to (\ref{no-other-vertex}), there exist
$\sigma_1,\sigma_2,\sigma_3\in \LS_{\hat\theta^*}$ such that
$\sigma_3(\hat\theta')<\hat
y'<\sigma_2(\hat\theta')<\hat{x}'<\sigma_1(\hat\theta')$ and except $\sigma_2$
there is no other line segment of $\cB$ between $\sigma_3$ and $\sigma_1$ over
the interval $(\hat\theta^*-\eta_*,\hat\theta^*)$. \medskip

Now, we show that $\CV_{\hat\theta^*}\neq \emptyset$ if
$\hat\theta^*<\infty$. Suppose this is not the case. Then for the above chosen
$\sigma_1,\sigma_2,\sigma_3\in LS_{\hat\theta^*}$, we have
$\sigma_1(\hat\theta^*)\neq \sigma_2(\hat\theta^*)\neq
\sigma_3(\hat\theta^*)$. Otherwise, the point
$(\hat\theta^*,\sigma_2(\hat\theta^*))$ is a critical vertex, which is a
contradiction to $\CV_{\hat\theta^*}=\emptyset$. Then for $\iota=1,2,3$, if
$\hat\theta^*\in\textrm{int}(\pi_1(\sigma_\iota))$, we let
$\sigma_{\iota+3}=\sigma_{\iota}$. Otherwise, there exists another
$\beta_\iota\in \LS_{\hat\theta^*}$ with
$\beta_\iota(\hat\theta^*)=\sigma_{\iota}(\hat\theta^*)$ and such that
$\hat\theta^*$ is the left endpoint of $\pi_1(\beta_\iota)$, and in this case we let
$\sigma_{\iota+3}=\beta_\iota$.  Then by (\ref{no-other-vertex}), except for
$\sigma_5$ there is no other line segment of $\cB$ between $\sigma_6$ and
$\sigma_4$ over the interval $(\hat\theta^*,\hat\theta^*+\eta_*)$.  Therefore, we can
extend $\Lambda^\pm$ continuously for $\hat\vartheta\in[\hat\theta^*-\eta_*/2,
  \hat\theta^*+\eta_*/2]$ with
$(\Lambda^-,\Lambda^+)\in\mathcal{CS}(\hat\theta,\hat x,\hat
y,\hat\theta^*+\eta_*/2)$ (see Figure \ref{Figure-extension} (a)).  This
contradicts the definition of $\hat\theta^*$.
\medskip

 Hence, we have $\CV_{\hat\theta^*}\neq \emptyset$. By Proposition~\ref{Step 3}, there
 exists exactly one critical vertex of $\tilde B$ on the fibre $\{\theta^*\}\times \T^1$,
 so there is exactly one critical vertex $(\hat\theta^*,\hat x^*)$ of $\cB$ in
 the line segment $\{\hat\theta^*\}\times[\sigma_3(\hat\theta^*),\sigma_3(\hat\theta^*)+1)$. We denote by
   $\alpha_1$ and $\alpha_2$ the two line segments of $\cB$ containing
   $(\hat\theta^*,\hat x^*)$.  If $\hat x^*\notin
   \{\sigma_\iota(\hat\theta^*)\ |\ \iota=1,2,3\}$, then we can extend $\Lambda^\pm$
   continuously on the interval $[\hat\theta^*-\eta_*/2, \hat\theta^*+\eta_*/2]$
   in the same way as in the situation with $\CV_{\hat\theta^*}=\emptyset$, which is
   a contradiction. Therefore, we have $\hat x^*\in
   \{\sigma_\iota(\hat\theta^*)\ |\ \iota=1,2,3\}$ and $(\hat \theta^*,\hat x^*)\in
   RCV_{\hat\theta^*}$ (note that $(\hat\theta^*,\hat x^*)$ cannot be a left
   critical vertex, since the line segments $\sigma_\iota$ with $\iota=1,2,3$
   extend to the left of $\hat\theta^*$). Furthermore, if $\hat x^*\neq
   \sigma_2(\hat \theta^*)$, then $\sigma_2\notin \{\alpha_1, \alpha_2\}$ and
   $\sigma_1(\hat\theta^*)\neq\sigma_2(\hat\theta^*)\neq\sigma_3(\hat\theta^*)$. Moreover,
   either $\hat x^*=\sigma_3(\hat \theta^*)$ or $\hat
   x^*=\sigma_1(\hat\theta^*)$. For the case $\hat x^*=\sigma_3(\hat \theta^*)$,
   let $\sigma_4,\sigma_5$ be as before.  Since $(\hat \theta^*,\hat x^*)\in
   \RCV_{\hat\theta^*}$, there is $\rho\in \LS_{\hat\theta^*}$ with
   $\rho(\hat\theta^*)<\sigma_3(\hat\theta^*)$, and except $\sigma_5$ there is
   no other line segment of $\cB$ between $\rho$ and $\sigma_4$ over the interval
   $(\hat\theta^*,\hat\theta^*+\eta_*)$ (see Figure \ref{Figure-extension} (b)).
   Thus, it is easy to see that we can extend $\Lambda^\pm$ continuously to the
   interval $[\hat\theta^*-\eta_*/2,\hat\theta^*+\eta_*/2]$ such that
   $(\Lambda^-,\Lambda^+)\in\mathcal{CS}(\hat\theta,\hat x,\hat
   y,\hat\theta^*+\eta_*/2)$, contradicting the definition of
   $\hat\theta^*$. The case where $\hat x^*=\sigma_1(\hat\theta^*)$ can be
   treated similarly.

  Therefore, we can assume that $\hat x^*=\sigma_2(\hat\theta^*)$ and
  $\sigma_2\in\{\alpha_1,\alpha_2\}$. Without loss of generality, we assume that
  $\alpha_2(\hat\vartheta)<\alpha_1(\hat\vartheta)$ for
  $\hat\vartheta\in(\hat\theta^*-\eta_*,\hat\theta^*)$.  Now we distinguish two
  cases.\medskip

  If $\sigma_2=\alpha_1$, then by the selection of $\sigma_1,\sigma_2,\sigma_3$,
  we have $\sigma_3=\alpha_2$ because there is no other line segment of $\cB$
  between $\alpha_2$ and $\alpha_1$ over $(\hat\theta^*-\eta_*,\hat\theta^*)$.
  By the definition of shadowing, for any $(\Lambda^-,\Lambda^+)\in
  \mathcal{CS}(\hat\theta,\hat x,\hat y,\hat\theta^*)$ we have that
  $\alpha_2(\hat\vartheta)<\Lambda^+(\hat\vartheta)<\alpha_1(\hat\vartheta)$ for all
  $\hat\vartheta\in (\hat\theta^*-\eta_*,\hat\theta^*)$.  Together with
  $\alpha_1(\hat\theta^*)=\alpha_2(\hat\theta^*)=\hat x^*$ (see Figure
  \ref{Figure-extension}~(c)), we obtain
  \begin{equation}
    \Lambda^+(\hat\theta^*) \ = \ \lim_{\hat\vartheta\rightarrow\hat\theta^{*}} {\Lambda}^+(\hat\vartheta)\ = \ \hat x^*.
   \end{equation}

   If $\sigma_2=\alpha_2$, then we have $\sigma_1=\alpha_1$ and similar to the
   above we conclude that
  \begin{equation}
    \Lambda^-(\hat\theta^*) \ = \ \lim_{\hat\vartheta\rightarrow\hat\theta^{*}} \Lambda^-(\hat\vartheta)=\hat x^* \ .
  \end{equation}
  This completes the proof.  \qed\medskip

  For $(\hat\theta,\hat x,\hat y)\in \cS$, if
  $\hat\theta^*=\hat\theta^*(\hat\theta,\hat x,\hat y)<\infty$, we say the point
  $(\hat\theta^*,\hat x^*)\in \RCV$ obtained in Lemma \ref{l.target-point} is
  the {\em target point} of $(\hat\theta,\hat x,\hat y)$.  Moreover, we denote
  by $\ell_*$ the maximal slope of a line segment in $\cB$ and note that since
  $\cB$ contains no vertical segments we have $\ell_*<\infty$.
Now we are ready to prove Proposition \ref{Step 4}.\medskip

\noindent\textbf{Proof of Proposition \ref{Step 4}.} Recall that by Definition \ref{d.mode-locking}, in
order to prove that a map in $\cF$ is mode-locked, we need to find a closed
topological annulus $\cA$ that is mapped strictly into itself by some iterate of
the map.

Suppose $\tilde f=(\tilde\omega,\tilde\varphi)\in\cF$ satisfies the assertions
of Proposition~\ref{Step 3}. Note that if we have $\hat\theta^*(\hat \theta,\hat
x,\hat y)=\infty$ for some $(\hat\theta,\hat x,\hat y)\in\cS$, then Lemma
\ref{l.target-point} implies that there exist closed curves $\gamma^\pm\subseteq
\tilde B^\pm$ in the same homotopy class. Let $\cA$ be the annulus that is bounded by
these curves and lies between $\gamma^+$ and $\gamma^-$ in the counterclockwise
direction. Then, since the lifts $\Gamma^+$ and $\Gamma^-$, which are the graphs of
continuous functions from $\R$ to $\R$, are mapped above, respectively below themselves
by the lift of $\tilde f^q$, the annulus $\cA$ is mapped strictly inside itself
by $\tilde f^q$.  Hence, we only need to address the case where
$\hat\theta^*(\hat \theta,\hat x,\hat y)<\infty$ for all $(\hat\theta,\hat
x,\hat y)\in\cS$.\medskip

Due to Proposition \ref{Step 3}, we have that $\epsilon=\min \{d(\vartheta',\vartheta'')
\mid (\vartheta',\xi')\neq (\vartheta'',\xi'')\in \pi(\CV)\}>0$.  Choose $(\hat
\theta,\hat x^-,\hat x^+)\in\cS$ such that for any $(\hat\vartheta,\hat\xi)\in \CV$ we
have $d(\hat\theta, \hat\vartheta)\geq \epsilon/2$.  We claim that for any $N\in\N$
there exist continuous curves
\begin{equation}\label{e.curves}
  \Gamma^\pm_N\ = \ \bigcup_{i=0}^N \Lambda_i^\pm\cup  \bigcup_{i=1}^N S_i^\pm
\end{equation}
such that for all $i=0\ld N$ we have
  \begin{enumerate}
  \item [(i)]
    $\pi_1(\Lambda_i^+)=\pi_1(\Lambda_i^-)=[\hat\theta_{i-1},\hat\theta_i]$ is a
    closed interval, denoted by $I_i$, with $\Lambda^\pm_{i|I_i}\subseteq \cB^\pm$, $|I_i|>\epsilon/3$,
    $\hat\theta_{-1}=\hat\theta$ and for $i\geq 0$ there exists
    $(\hat\theta_i^*,\hat x_i^*)\in\RCV$ such that
    $\hat\theta_i=\hat\theta_i^*-\epsilon/10$;
  \item [(ii)]
    $(\Lambda_i^-,\Lambda^+_i)\in\cC\cS(\hat\theta_{i-1}+\epsilon/5,\Lambda^-_i(\hat\theta_{i-1}+\epsilon/5),
    \Lambda^+_i(\hat\theta_{i-1}+\epsilon/5),\hat\theta_i)$ and $(\hat\theta_i^*, \hat
    x_i^*)$ is the target point of $(\hat\theta_{i}-\epsilon/10,
    \Lambda_i^-(\hat\theta_{i}-\epsilon/10),\Lambda_i^+(\hat\theta_{i}-\epsilon/10))\in
    \cS$;
  \item [(iii)] either $S_{i+1}^+=\emptyset$,
    $\Lambda_i^+(\hat\theta_i)=\Lambda_{i+1}^+(\hat\theta_i)$ and $S_{i+1}^-$ is a
    vertical segment over $\hat\theta_i$ connecting $\Lambda^-_i$ and
    $\Lambda^-_{i+1}$ with $\Lambda^-_i(\hat\theta_i)<\Lambda^-_{i+1}(\hat\theta_{i})$ , or
    $S_{i+1}^-=\emptyset$, $\Lambda_i^-(\hat\theta_i)=\Lambda_{i+1}^-(\hat\theta_i)$
    and $S_{i+1}^+$ is a vertical segment over $\hat\theta_i$ connecting
    $\Lambda_i^+$ and $\Lambda_{i+1}^+$ with
    $\Lambda_i^+(\hat\theta_i)>\Lambda_{i+1}^+(\hat\theta_i)$.
    \end{enumerate}

    In order to start the inductive construction, we let
    $\hat\theta_{-1}=\hat\theta$ and apply Lemma \ref{l.target-point} to obtain
    $(\bar\Lambda_0^-,\bar\Lambda_0^+)\in\mathcal{CS}(\hat\theta,\hat x^-,\hat
    x^+,\hat \theta_0^*)$, where $(\hat \theta_0^*, \hat x_0^*)$ is the target
    point of $(\hat\theta,\hat x^-,\hat x^+)$.  Then we let
    $I_0=[\hat\theta,\hat\theta_0^*-\epsilon/10]$ and
    $\Lambda_0^\pm=\bar\Lambda^\pm_{0|I_0}$. It is easy to check that (i) and (ii)
    hold, whereas (iii) is still void.\medskip

    Now, assume that we have already constructed $\Lambda_i^\pm$ for
    $i=0,1,\ldots,\iota$ and $S_i^{\pm}$ for $i=1,\ldots, \iota$. Since
    $\hat\theta_\iota=\hat\theta^*_\iota-\epsilon/10$ by (i) and due to the choice
    of $\epsilon$, we know that there is no critical vertex in
    $(\hat\theta_\iota-\epsilon/5,\hat\theta_\iota^*)\times \R$. Further,
    $(\hat\theta_\iota^*,\hat x_\iota^*)$ is the target point of
    $(\hat\theta_{\iota}-\epsilon/10,
    \Lambda_\iota^-(\hat\theta_{\iota}-\epsilon/10),
    \Lambda_\iota^+(\hat\theta_{\iota}-\epsilon/10))$. Hence, by Lemma
    \ref{l.target-point} we can continuously extend $\Lambda_\iota^\pm$ to
    continuous curves $\bar\Lambda_\iota^\pm$ defined on the larger interval
    $[\hat\theta_{\iota-1},\hat\theta_\iota^*]$. (In fact, the $\Lambda^\pm_i$ are
    defined as the restrictions of such longer curves in the inductive
    construction, see below.) There are two possibilities
\begin{itemize}
\item [(1)] $\bar\Lambda_\iota^-(\hat\theta_\iota^*)=\hat x_\iota^*$;
\item [(2)] $\bar\Lambda_\iota^+(\hat\theta_\iota^*)=\hat x_\iota^*$.
\end{itemize}

 \begin{figure}[ht!]
 \epsfig{file=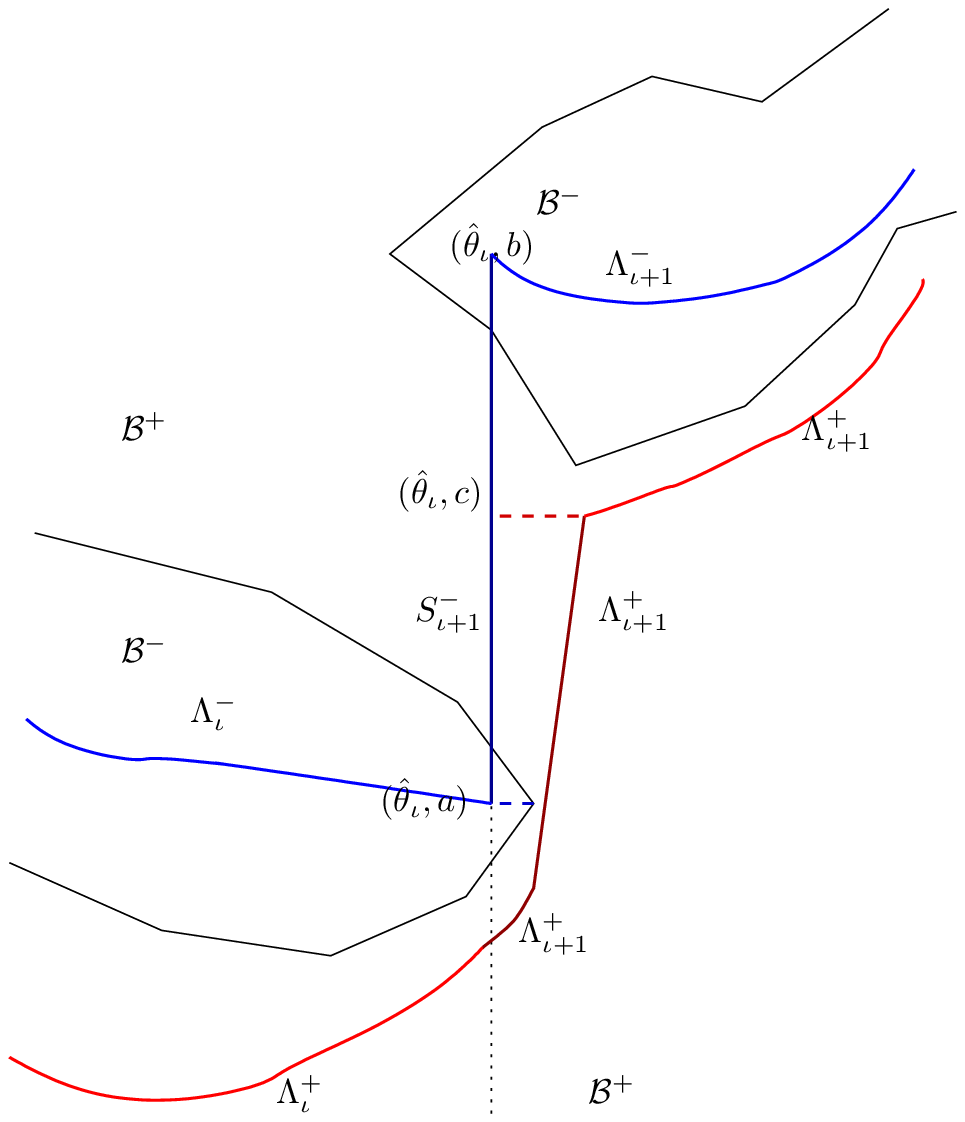,
 clip=,width=0.45\linewidth} \hspace{2eM}
 \epsfig{file=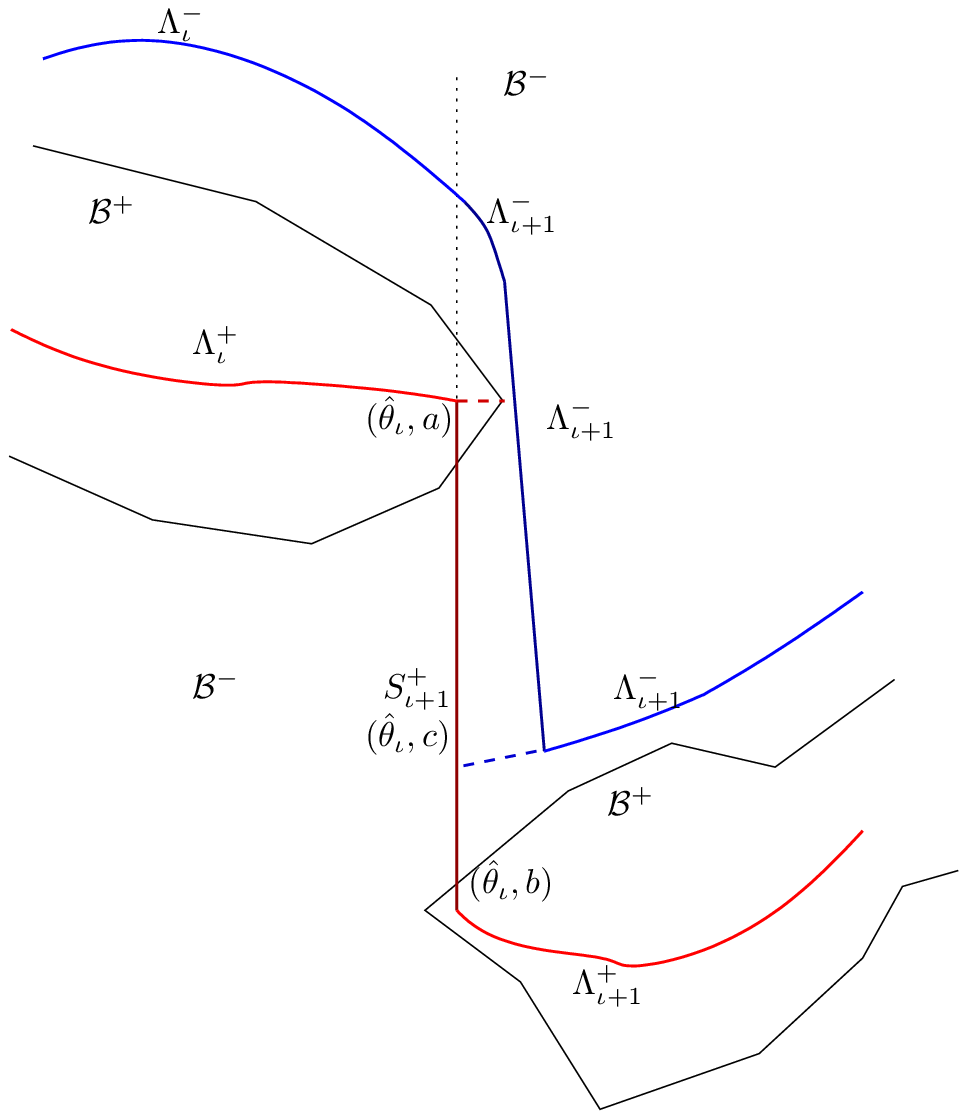, clip=,width=0.44\linewidth}

\vspace{1ex}
\hspace{6eM} Case (1)\hspace{18eM} Case (2)
 \vspace{1ex}

\caption{The strategy of the construction. \label{Figure-construction} }
\end{figure}

\textit{Case (1):} If $\bar\Lambda_\iota^-(\hat\theta_\iota^*)=\hat x_\iota^*$,
we let $S_{\iota+1}^+=\emptyset$ and choose $S_{\iota+1}^-$ to be a vertical segment with
endpoints $(\hat\theta_\iota,a)$, $(\hat\theta_\iota,b)$ where
$a=\Lambda_\iota^-(\hat\theta_\iota)$ and $b>a$ is chosen such
that $(\hat\theta_\iota,a)$ and $(\hat\theta_\iota,b)$ belong to two
successive connected components of
$\cB_{\hat\theta_\iota}^-=(\{\hat\theta_\iota\}\times \R)\cap \cB^-$. Since
there is no critical vertex with first coordinate $\hat\theta_\iota$, there
exists a point $(\hat\theta_\iota,c)\in \cB^+$ in the shadow of
$(\hat\theta_\iota,b)$.  Then, by Lemma \ref{l.target-point}, there exists
a unique target point $(\hat \theta_{\iota+1}^*, \hat x_{\iota+1}^*)$ of
$(\hat\theta_\iota,b,c)\in\cS$ and there is
$(\bar\Lambda_{\iota+1}^-,\bar\Lambda_{\iota+1}^+)\in\mathcal{CS}
(\hat\theta_\iota,b,c,\hat\theta_{\iota+1}^*)$. By
the choice of $b$, we have that
$\hat\theta^*_{\iota+1}\neq\hat\theta^*_\iota$ and hence
$\hat\theta^*_{\iota+1}>\hat\theta^*_\iota$.

 Since there is no other critical vertex in
 $B_{\epsilon/2}(\theta_\iota^*)\times\T^1$, there exists $0<\delta\leq
 \epsilon/10$ such that the line segment $\lambda$ connecting the points
 $(\hat\theta_\iota^*, \bar\Lambda_\iota^+(\hat\theta_\iota^*))$ and $(\hat
 \theta_\iota^*+\delta, \bar\Lambda_{\iota+1}^+(\hat\theta_\iota^*+\delta))$ is
 contained in $\cB^+$ (see Figure~\ref{Figure-construction}(Case (1))). Now we let
 $\hat\theta_{\iota+1}=\hat\theta_{\iota+1}^*-\epsilon/10$,
 $I_{\iota+1}=[\hat\theta_\iota,\hat\theta_{\iota+1}]$ and set
 $\Lambda^-_{\iota+1}=\bar\Lambda^-_{\iota+1|I_{\iota+1}}$. Further, we define
 $\Lambda^+_{\iota+1}$ as
 \[
\Lambda^+_{\iota+1} \ = \ \bar\Lambda^+_{\iota|[\hat\theta_\iota,\hat\theta_\iota^*]}\cup
\lambda \cup \bar\Lambda^+_{\iota+1|[\hat\theta_\iota^*+\delta,\hat\theta_{\iota+1}]} \ .
 \]
 It is now easy to check that $\Lambda^\pm_{\iota+1}$ and $S^\pm_{\iota+1}$
 satisfy the inductive assumptions (i)--(iii). Moreover, in Case (2) the
 construction works in completely symmetric way. \medskip

Hence, we have that curves $\Gamma^\pm_N$ with properties (i)--(iii) exist for
all $N\in\N$. Now, suppose $M=\#\pi(\RCV)$. Then there have to exist $i< j\in
\{0\ld M\}$ such that $\pi(\hat\theta_{i}^*,\hat x_{i}^*)=\pi(\hat\theta_{j}^*,
\hat x_{j}^*)$ and for any other pair $(i',j')$ with $i'\neq
j'\in\{i,i+1,\ldots, j\}$, $\pi(\hat\theta_{i'}^*,\hat
x_{i'}^*)\neq\pi(\hat\theta_{j'}^*, \hat x_{j'}^*)$. Without loss of generality,
we assume $i=0$.  Then, we have that $\hat\theta^*_j=\hat\theta^*_0+k$ and $\hat
x^*_j=\hat x^*_0+m$ for some integers $k,m$, so the target points of
$\Lambda_0^\pm$ and $\Lambda_j^\pm$ only differ by the integer vector
$(k,m)$. This means that the projections of these curves to $\ntorus$ have to be
contained in the same connected components of $\tilde B^\pm$. Since the sets
$\cB^\pm$ are open, we can thus modify $\Lambda_j^\pm$ in such a way that for
some $\tilde\delta>0$ we have $\Lambda_j^\pm(\hat\vartheta+j)=\Lambda^\pm_0(\hat\vartheta)+k$
for all $\hat\vartheta\in[\hat\theta_0-\tilde\delta,\hat\theta_0]$. As a consequence, the
curves
\[
\bar\Gamma^\pm \ = \ \bigcup_{i=1}^N \Lambda_i^\pm \cup \bigcup_{i=1}^{N}S^\pm_i
\]
project to closed curves $\gamma^\pm$ in the torus. Note that due to the
modification of $\Lambda_j^\pm$, the projections of these curves end exactly in
the starting points of the projection of $\Lambda_1^\pm$ or $S_1^\pm$,
respectively. We denote by $\Gamma^\pm$ the periodic extensions of the curves
$\bar\Gamma^\pm$ to all of $\R$, that is, the $\Gamma^\pm$ are lifts of
$\gamma^\pm$ to $\R^2$. \medskip

Finally, we want to see that the existence of the curves $\Gamma^\pm$
constructed above implies that for any $\omega'$ in a right neighbourhood of
$\tilde\omega$ the system $\bar f=(\omega',\tilde\varphi)$ is mode-locked. This
follows if we can show that a suitable lift $\bar F^q$ of $\bar f^q$ maps the
`infinite strip'
\[
\wh \cA \ = \ \{(\hat\vartheta,\hat\xi)\in\R^2\mid \Gamma^+(\hat\vartheta)\leq \hat\xi \leq \Gamma^-(\hat\vartheta)\}
\]
strictly inside itself. This forces the fibred rotation number of $\bar f$ to be
equal to $\frac{m}{kq}$, and since this situation is persistent under
perturbations of $\tilde\varphi$ it implies that $\bar f$ is mode-locked
(compare \cite{bjerkloev/jaeger:2009}).

In order to do so, we choose $\tilde F^q$ to be a suitable lift of $\tilde f^q$
to $\R^2$ which has rotation number zero in the base, and for $\omega'$ slightly
to the right of $\tilde \omega$ we choose $\bar F^q$ as the lift of $\bar f^q$
whose rotation number is close to zero. (Note here that we choose $\tilde F^q$
directly as a lift of $\tilde f^q$, and it is not necessarily the iterate of a
lift $\tilde F$ of $\tilde f$, similar to $\bar F^q$.) Note that for $\tilde
F^q$ itself, the vertical parts of the curves $\Gamma^\pm$ contain fixed points
(since they intersect both $\cB^-$ and $\cB^+$), so that $\cA$ cannot be mapped
strictly inside itself. However, if $\omega'$ is slightly larger than
$\tilde\omega$ then $\bar F^q$ shifts these vertical segments slightly to the
right, so that they become disjoint from their images. It remains to show that
the whole curves $\Gamma^\pm$ are also disjoint from their images under $\bar
F^q$. The fact that $\Gamma^-$ is mapped below and $\Gamma^+$ is mapped above
themselves then follows easily from the fact that this is true for the
non-vertical parts of these curves under the map $\tilde F^q$.\medskip

We give the argument for $\Gamma^-$, while the proof for $\Gamma^+$ is
analogous. Suppose that $\omega'\in(\tilde\omega,\tilde\omega+s)$, where $s>0$
is small (to be specified below).  In order to prove that $\bar F^q(\Gamma^-)
\cap \Gamma^-=\emptyset$, we need to show that
\begin{enumerate}
\item[(1)] $\bar F^q(\Lambda^-_i)\cap \Lambda^-_j=\emptyset$ for all
  $i,j\in\{1\ld N\}$.
\item[(2)] $\bar F^q(\Lambda^-_i)\cap S^-_j=\emptyset$ for all $i,j\in\{1\ld N\}$.
\item[(3)] $\bar F^q(S^-_j)\cap \Lambda^-_i=\emptyset$ for all $i,j\in\{1\ld
  N\}$.
\end{enumerate}
\begin{enumerate}
\item[(1)] If $|i-j|\geq 1$, then $\bar F^q(\Lambda^-_i)$ and $\Lambda^-_j$ have
  disjoint projections to the first coordinate, since this is true for
  $\Lambda^-_i$ and $\Lambda^-_j$ and $\bar F^q$ only shifts horizontally by
  $q(\omega'-\tilde\omega)\leq qs$, which is arbitrarily small. This is similar
  for $i=j+1$, since $\bar F^q$ shifts horizontally to the right.

  If $i=j$, then $\tilde F^q(\Lambda^-_i)$ and $\Lambda^-_i$ are disjoint by
  construction ($\Lambda^-_i\ssq\cB^-$) and therefore have positive distance to
  each other. Hence, if $s$ is sufficiently small, this remains true if $\tilde
  F^q$ is replaced by $\bar F^q$.

  Finally, if $i=j-1$, then the projections of $\bar F^q(\Lambda^-_i)$ and
  $\Lambda^-_j$ have a small interval
  $[\hat\theta_i,\hat\theta_i+q(\omega'-\tilde\omega)]\subseteq[\hat\theta_i,\hat\theta_i+qs]$
  of overlap. We have to consider two situations: If $S^-_j=\emptyset$, then
  $\Lambda^-_i\cup\Lambda^-_j$ is a connected curve that is mapped strictly
  below itself by $\tilde F^q$. As above, this remains true after perturbed to
  $\bar F^q$ if $s$ is sufficiently small. If $S^-_j\neq \emptyset$, then by
  construction the endpoints $(\hat\theta_i,\Lambda_i^-(\hat\theta_i))$ and
  $(\hat\theta_i,\Lambda_{j}^-(\hat\theta_i))$ of this segment are separated by
  some interval that belongs to $\cB^+$. Hence, $\tilde
  F^q(\hat\theta_i,\Lambda_{j}^-(\hat\theta_i))$ lies strictly above
  $(\hat\theta_i,\Lambda_i^-(\hat\theta_i))$. By continuity, this implies that
  the overlapping parts of the curves $\bar F^q(\Lambda^-_i)$ and $\Lambda^-_j$
  are still disjoint if $s$ is sufficiently small.
\item[(2)] Similar to before, since $\bar F^q$ only shifts horizontally by at
  most $qs$, we only have to consider the cases $i=j-1$ and $i=j$. If $i=j-1$,
  then $\Lambda_i^-$ connects to $S^-_j$ in its right endpoint, which is the
  lower endpoint of $S^-_j$. As $\Lambda^-_i$ is mapped below itself by $\tilde
  F^q$ its image will be disjoint from $S^-_j$, and the same is still true after
  perturbation if $s$ is small.

  If $i=j$, then the upper endpoint of $S^-_j$ equals the left endpoint of
  $\Lambda^-_i$. Since $\bar F^q$ shifts slightly to the right, this means that
  the projection $\hat\theta_{j-1}$ of $S^-_j$ is not contained in the
  projection of $\bar F^q(\Lambda^-_i)$, so the two sets are disjoint.
\item[(3)] Again, we only have to consider $i=j-1$ and $i=j$. In the first case,
  $S^-_j$ is mapped to the right by $\bar F^q$ and therefore has disjoint
  projection from $\Lambda_i^-$. In the second case, $S_j^-$ lies below the left
  endpoint of $\Lambda_j^-$. Since the upper endpoint of $S_j^-$ is mapped
  downwards by $\tilde F^q$, the image $\tilde F^q(S_j^-)$ has positive distance
  to $\Lambda_i^-$, and as before this remains true after perturbation if $s$ is
  small.\qed
\end{enumerate}
\medskip

\section{Proof of Corollary~\ref{c.main-homeomorphisms} and
  Theorem~\ref{t.main-families}}
\label{FurtherProofs}

Corollary \ref{c.main-homeomorphisms} is a direct consequence of Theorem
\ref{t.main-homeomorphisms} and the following elementary lemma on products of
Baire spaces.  

\begin{Lemma}
  Suppose $A, B$ are Baire spaces, B is second-countable, and $U\subseteq
  A\times B$ is a residual (resp.\ open and dense) set. Then there exists a
  residual set $V\subseteq A$, such that for any $v\in V$ there exists a
  residual (resp.\ open and dense) set $U_v\subseteq B$ such that $\{v\}\times
  U_v\subseteq U$.
\end{Lemma}


%

 In the same way, Theorem \ref{t.main-families}(b) follows from Theorem
 \ref{t.main-families}(a). Hence, we only need to prove Theorem
 \ref{t.main-families}(a).

\proof[\bf{Proof of Theorem \ref{t.main-families}(a)}\ :\ ]
We fix an arbitrary dense sequence $\{\tau_n\}_{n\in\N}\subseteq \T^1$ and let
\[\cP_n \ =\ \left\{(f_\tau)_{\tau\in\T^1}\in \cP\ |\ f_{\tau_n} \textrm{ is mode-locked}\ \right\}.\]
If all $\cP_n$ are open and dense in $\cP$, then $\cP^{DS}=\cap_n \cP_n$ is
residual and satisfies the assertion. Here, it suffices to show that $\cP_n$ is
dense, openness being obvious.

To that end, we fix $(f_\tau)_{\tau\in\T^1}\in \cP$, and $\varepsilon>0$.
Choose $\hat f\in B_{\varepsilon/3}(f_{\tau_n})$ so that $\hat f$ is
mode-locked. We write $\hat f=(\hat \omega, \hat \varphi), f_\tau=(\omega,
\varphi_\tau)$. Then we fix an interval $I=[a,b]$ with
$\tau_n\in\textrm{int}(I)$, such that $f_\tau\in B_{\varepsilon/3}(f_{\tau_n}),
\ \forall \tau\in I$.  We define
\[\tilde f_\tau=\left\{\begin{array}{ll}(\hat \omega, \varphi_\tau)& \textrm{if}\  \tau\notin I,\\
    (\hat \omega, \frac{\tau_n-\tau}{\tau_n-a}\varphi_a+
    \frac{\tau-a}{\tau_n-a}\hat \varphi )& \textrm{if}\  \tau\in [a,\tau_n],\\
    (\hat \omega, \frac{b-\tau}{b-\tau_n}\hat \varphi+ \frac{\tau-\tau_n}{b-\tau_n}\varphi_b)
    & \textrm{if} \ \tau\in [\tau_n, b],
\end{array}\right.\]
where the convex combinations should be interpreted as discussed in Section
\ref{Preliminaries}.  Then $(\tilde f_\tau)_{\tau\in\T^1}$ is
$\varepsilon$-close to $(f_\tau)_{\tau\in\T^1}$ and $\tilde f_{\tau_n}=\hat f$
is mode-locked. This completes the proof. \qed

\end{document}